\documentclass[12pt]{article}
\usepackage{amscd, amsfonts, amsmath, amssymb, amsthm, color}
\usepackage[all]{xy}
\usepackage{pinlabel}
\usepackage{stmaryrd}
\usepackage{graphicx}
\usepackage{hyperref}
\newtheorem{thm}{Theorem}[section]
\newtheorem{prp}[thm]{Proposition}
\newtheorem{cor}[thm]{Corollary}
\newtheorem{lma}[thm]{Lemma}

\theoremstyle{definition}
\newtheorem{dfn}[thm]{Definition}
\newtheorem{rmk}[thm]{Remark}

\numberwithin{equation}{section}

\newcommand{\R}{{\mathbb{R}}}
\newcommand{\C}{{\mathbb{C}}}
\newcommand{\Q}{{\mathbb{Q}}}
\newcommand{\Z}{{\mathbb{Z}}}

\newcommand{\bbH}{{\mathbb{H}}}

\newcommand{\Cc}{{\mathcal{C}}}

\newcommand{\Ee}{{\mathcal{E}}}

\newcommand{\Ss}{{\mathcal{S}}}

\newcommand{\Xx}{{\mathcal{X}}}
\newcommand{\cH}{{H}}

\newcommand{\M}{{\mathcal{M}}}

\newcommand{\la}{\langle}

\newcommand{\ra}{\rangle}

\newcommand{\pa}{\partial}

\newcommand{\ju}{{\operatorname{jp}}}

\newcommand{\diag}{{\operatorname{diag}}}
\newcommand{\BV}{\operatorname{BV}}

\newcommand{\wt}{\widetilde}
\newcommand{\wh}{\widehat}

\newcommand{\ul}{\underline}
\newcommand{\p}{\partial}

\newcommand{\eps}{\varepsilon}

\newcommand{\oX}{\overline{X}}

\newcommand{\bU}{{\mathbf U}}
\newcommand{\bM}{{\mathbf M}}
\newcommand{\st}{\operatorname{st}}

\newcommand{\cyc}{{\rm cyc}}

\newcommand{\Id}{\mathrm{id}}%Changed from Id to id (TE)

\newcommand{\Ho}{\mathrm{Ho}}
\newcommand{\bal}{{\mathrm{bal}}}

\newcommand{\MM}{\mathcal{M}}

\newcommand{\CC}{\mathcal{C}}

\newcommand{\XX}{\mathcal{X}}

\newcommand{\lb}{\left(}
\newcommand{\rb}{\right)}

\newcommand{\motimes}{\mathop{\otimes}}
\newcommand{\mstar}{\mathop{\star}\limits}
\newcommand{\llb}{\llbracket}%{\llparenthesis}
\newcommand{\rrb}{\rrbracket}%{\rrparenthesis}

\newcommand{\carrl}{\circlearrowleft}
\newcommand{\uda}{{\uparrow\downarrow}}
\newcommand{\dua}{{\downarrow\uparrow}}
 \newcommand{\uc}{{\underline{c}}}
  
   \newcommand{\ua}{{\underline{a}}}
    
        \newcommand{\ux}{{\underline{x}}}
\title
{Symplectic homology product via Legendrian surgery}
\author{Fr{\'e}d{\'e}ric Bourgeois\footnote{Partially supported by the Fonds National de la Recherche Scientifique (Belgium)}
\\Universit\'e Libre de Bruxelles \and Tobias Ekholm\footnote{Partially supported by the G{\"o}ran Gustafsson Foundation for Research in Natural Sciences and Medicine}\\Uppsala University \and
Yakov Eliashberg\footnote{Partially supported by NSF grants DMS-0707103 and DMS 0244663}
\\ Stanford University}

\begin{document}
\maketitle

\begin{abstract}
This research announcement continues the study of the symplectic homology of Weinstein manifolds undertaken in \cite{BEE1} where the symplectic homology, as a vector space, was expressed in terms of the Legendrian homology algebra of the attaching spheres of critical handles. Here we express the product and $\BV$-operator of symplectic homology in that context.      
\end{abstract}

%\tableofcontents

\section{Introduction}
Weinstein manifolds are symplectic counterparts of affine (Stein) complex manifolds. They are exact symplectic manifolds which admit symplectic handle decompositions with isotropic core disks. In particular the index of any handle is at most half of the dimension of the manifold. We refer the reader to \cite{BEE1,We,EG} for precise definitions.
 Middle dimensional handles are called {\em critical} and handles below middle dimension {\em subcritical}. Given a Weinstein manifold $X$, its symplectic homology $S\bbH(X)$ is a symplectic invariant which was first defined by A.~Floer and H.~Hofer in \cite{FlHo}. It vanishes if $X$ is subcritical, see \cite{Cie}. A general Weinstein manifold $X$ is obtained by attaching critical handles to a  
subcritical manifold $X_0$ along a collection $\Lambda$ of Legendrian spheres in the contact manifold $Y_0$ which is the ideal boundary of $X_0$. In \cite{BEE1}, $S\bbH(X)$ was expressed in terms of the Legendrian homology algebra of $\Lambda\subset Y_0$.

Symplectic homology comes with TQFT-like operations which on the homology level are all expressed in terms of a (pair-of-pants) product and a BV-operator, see \cite{Seidel-biased, Getzler}. In this paper we express the product and the BV-operator in terms of the Legendrian homology algebra of $\Lambda$ and  some additional operations on it. The BV-operator is determined by the algebra itself,  whereas the expression of the product involves a duality operation which can be identified with the first term in the extension of Legendrian homology algebra to Legendrian rational symplectic field theory (SFT). 

The paper is organized as follows. In Section \ref{sec:alg-prelim} we introduce algebraic constructions which we later apply to Legendrian homology algebras.  In Section \ref{sec:leg-algebra} we express the Legendrian homology algebra results from \cite{BEE1} in the terminology of Section \ref{sec:alg-prelim} and also include a discussion on the construction of a version of rational SFT in the Legendrian setting  (completion of this construction is work in progress). In Section \ref{sec:BV-Leg} we present the symplectic homology product and the BV-operator in terms of Legendrian homology algebra and a duality operator and in Section \ref{sec:examples} we use the result to compute the symplectic homology, with product and BV-operator, of cotangent bundles of spheres.  While proofs are mostly omitted, we note that the  proofs of the algebraic results from  Section \ref{sec:alg-prelim}  are  fairly straightforward.
\medskip

The authors thank M.~Abouzaid, S.~Ganatra, E.~Getzler, L.~Ng, and P.~Seidel for inspiring discussions, and Y.~Lekili for correcting a sign error in \eqref{eq:dM}.

\section{Algebraic preliminaries}\label{sec:alg-prelim}
In this section we describe a number of constructions in a purely algebraic setting. The constructions will be applied in geometrically relevant situations in Sections   \ref{sec:leg-algebra} and  \ref{sec:BV-Leg}.
\subsection{A tensor algebra and associated  objects}\label{sec:MoM}
We associate to a differential graded algebra $A$ over a field $K$ a pair of dual $A$-modules and a corresponding tensor algebra. We also discuss cyclic versions of these objects, which are vector spaces over $K$.
\subsubsection{The algebra $A$}
Let $R$ be an algebra over a field $K$ of characteristic $0$ generated by idempotents $1_1,\dots, 1_m$, where
$1_k1_l=\delta_{kl}1_k$, where $\delta_{kl}$ is the Kronecker delta. 
We define a  graded module  $V$ over  $R$, which is generated by a countable set $\Cc$ decomposed as a disjoint union $\Cc=\bigcup_{1\le i,j\le m} \Cc_{ij}$. The degree of an element $a\in V$ will be denoted by $|a|$.
The decomposition of $\Cc$ induces a decomposition 
$V=\bigoplus\limits_{1\leq i,j\leq m} V^{(ij)}$, where $V^{(ij)}$ is generated by $\Cc_{ij}$. When applied on the left, the element $1_j$, $j=1,\dots, m$, acts as the identity  on $V^{(ij)}$ and as the $0$-map on $V^{(il)}$, $l\ne j$. Similarly, when applied on the right, the element $1_j$  acts as the  identity on $V^{(ji)}$ and as $0$ on $V^{(li)}$, $l\ne j$.

We define an extension $\wh V$ of the module $V$ as follows.  The module $\wh V$ 
is generated by $\wh\CC\cup \XX$, where $\wh\CC=\{\wh c\colon c\in\CC\}$ and 
$\Xx=\{x_1,\dots,x_m\}$, where $|\wh c|=|c|+1$, for $c\in\CC$ and $|x_j|=0$, $j=1,\dots,m$.
We define 
\[
\wh V :=K\la \Xx\ra \oplus V[1],
\]
where $K\la \Xx\ra$ denotes the $m$-dimensional module freely generated by the set $\Xx$ and where $V[1]$ denotes $V$ with grading shifted up by $1$.  The module $\wh V$ inherits a  decomposition $\wh V=\bigoplus\limits_{1\leq i,j\leq m} \wh V^{(ij)}$  from $V$, where 
\[
\wh V^{(ij)}=
\begin{cases}
V^{(ij)}[	1] &\text{if } i\ne j,\\
Kx_j\oplus V^{(jj)}[1] &\text{if }i=j.
\end{cases}
\] 
 Here the idempotents act as before on $V\approx V[1]$ as a subset of $\wh V$, and as follows on the new generators: $1_jx_i=\delta_{ji}x_j$ and $x_i1_j=\delta_{ij}x_j$, $1\le i,j\le m$.

We will use the following notation. If $a\in V$ then $\wh a$ denotes the element in $V[1]\subset\wh V$ that corresponds to $a$. Furthermore, if $F$ and $G$ are either $V$ or $\wh V$ and if  
$(u, v)$ is an ordered pair of elements in $F^{(ij)}\oplus G^{(kl)}$ then $(u,v)$ is {\em composable} if $j=k$.

Let 
\begin{equation}\label{eq:defA}
A :=R\oplus V\oplus V^{{\mathop{\otimes}\limits_R}^2}\oplus\dots
\end{equation}
be the tensor algebra over $R$ generated by $V$ with its natural grading.  Then $A$ is generated additively by composable monomials $w=c_1\dots c_k$, $c_j\in V$, $j=1,\dots,k$, where we say that $w$ is composable provided $(c_j,c_{j+1})$ is composable for $1\le j\le k-1$. (Here, as often in the following, we will suppress the tensor sign in the notation for tensor products). The decomposition of $V$ induces a corresponding decomposition of $A$, $A=\bigoplus_{1\le i,j\le m} A^{(ij)}$, where $A^{(ij)}$ is generated by composable monomials $c_1\dots c_k$ with $c_1\in V^{(is)}$ for some $s\in\{1,\dots, m\}$ and $c_k\in V^{(tj)}$ for some $t\in\{1,\dots,m\}$.
 
We will assume that the algebra $A$ has a differential $d\colon A\to A$ of degree $-1$ which leaves the subspaces $A^{(ij)}$, $i,j\in\{1,\dots,m\}$ invariant (i.e., $d^2=0$,  $d$ satisfies the Leibniz rule, and $d(A^{(ij)})\subset A^{(ij)}$).   
%For $a\in V$ we  write $da=\sum_{k=0}^{N(a)} d_ka, $ where $d_ka\in V^{{\mathop{\otimes}\limits_R}^k}$.
We will usually write $H(A)$ for the homology algebra $H_*(A,d)$.
 
\subsubsection{The $A$-module $M$ and the $K$-complex  $M^\cyc$}\label{sec:MMcyc}
 Define 
\begin{equation}\label{eq:defM}
M= M(A):=A\mathop{\otimes}\limits_R\wh V\mathop{\otimes}\limits_R A.
\end{equation}
Then $M$ is a graded left-right module over the differential graded algebra $A$ which is additively generated by composable monomials  $w=c_1\dots c_k u f_1\dots f_l$, where $c_i\in \Cc$, $i=1,\dots,k$, $u\in \wh \Cc\cup\Xx$, and $f_j\in \Cc$, $j=1,\dots,l$. We say that a composable monomial $w$ is {\em cyclically composable} if the pair $(f_l,c_1)$ is composable. Let $M^{\diag}\subset M$ be the submodule generated (over $K$) by cyclically composable monomials and define $M^{\cyc}$ as the quotient space of $M^{\diag}$  in which monomials which differ by graded cyclic permutation are identified.  If $w\in M^{\diag}$ then we write $\llb w\rrb$ for its equivalence class in  $M^{\cyc}$, i.e.~if $\pi\colon M^{\diag}\to M^\cyc$ denotes the natural projection then $\llb w\rrb:=\pi(w)\in M^\cyc$.

The differential $d\colon A\to A$ induces a  differential  $d_{M}\colon M\to M$ as follows.
If $w_1 u\, w_2\in  M$, $u\in\wh V$, and $w_1,w_2\in A$, then
\begin{equation}\label{eq:S1}
d_ {M}(w_1u\,w_2):=
(dw_1) u\,w_2+ (-1)^{|w_1|} w_1 (d_ {M} u) w_2+(-1)^{|w_1u|} w_1u (dw_2).
\end{equation}
Here
\begin{equation}\label{eq:dM}
d_M u=
\begin{cases}
ax_i- x_ja - S(da) &\text{if }u=\wh a,\,\,a\in V^{(ij)},\\
0 &\text{if }u=x_j,\;j\in\{1,\dots,m\},
\end{cases}
\end{equation}
where the $R$-module homomorphism $S:A\to M$ is defined by
\begin{align}\notag
S(a_1\dots a_k)&=\wh a_1a_2a_3\dots a_k\,\,+\,\,(-1)^{|a_1|} a_1\wh a_2a_3\dots a_k\\\label{eq:Sop}
&+\,\,\dots\,\,+\,\, (-1)^{|a_1\dots a_{k-1}|} a_1a_2\dots a_{k-1}\wh a_k,
\end{align}
and $S(1_j) = 0$, for $j= 1, \ldots, m$.
\begin{lma}\label{lm:dM2}
The map $d_M\colon M\to M$ is a differential (i.e.~$d_{M}^{\;\,2}=0$) which leaves $M^{\diag}$ invariant and descends to a differential $M^{\cyc}\to M^{\cyc}$, still denoted $d_{M}$, in such a way that the natural projection 
$\pi\colon (M^{\diag},d_ {M})\to (M^\cyc,d_{M})$ is a chain map.
\end{lma}
  
We write $H(M)$ and $H(M^\cyc)$, respectively,  for the homologies $H_*(M,d_{M})$ and  
$H_*( M^\cyc,d_ {M})$. Note that $H(M)$ is a left-right module over $H(A)$, while $H(M^\cyc)$ is just a $K$-vector space. 

\begin{rmk}\label{rem:satellites}
We can write 
$$
M^\cyc=\lb\bigoplus_{j=1}^m Ax_j\rb\oplus A\motimes\limits_R V[1] :={'}\!M^\cyc \oplus {''}\!M^\cyc
$$ 
(i.e.~${'}\!M^{\cyc}$ is the submodule of $M^{\cyc}$ generated by $\XX$ and ${''}\!M^{\cyc}$ its complement).
Let $p\colon M^\cyc\to {''}\!M^\cyc$ be the projection and let ${''}\!d_{M^\cyc}:=p\circ (d_{M^\cyc})|_{{''}\!M^\cyc}$. The second summand ${''}\!M^\cyc$
 can be viewed as a non-commutative version of  the ``Lie algebra of  vector fields''  on the space $V$.  
Indeed, we can identify the cyclic monomial $\llb w\wh u\rrb\in M^{\cyc}$, $w\in A$, $u\in V$, 
with the vector field $w\,\frac{\p}{\p u}$. The differential  ${''}\!d_{M^\cyc}$
 can be naturally viewed as a vector field on $V$, and hence an element of  ${''}\!M^\cyc$, given
 by $\sum_{c \in \Cc} dc \,\frac{\p}{\p c}$. 
   For  any $X\in {''}\!M^\cyc$  
   its differential ${''}\!d_MX$ is just the Lie bracket $ [{''}\!d_M,X]$ of vector fields, or equivalently the Lie derivative $L_{{''}\!d_M}X$.
   
This interpretation corresponds to the ``satellite philosophy'' from \cite{SFT}. Another interpretation  was pointed out to the authors by  M.~Abouzaid: $H(M^\cyc)$ is isomorphic to the Hochschild homology $HH_*(A,A)$.\end{rmk}

\subsubsection{The extended  tensor algebra $U$ and the space $U^\bal$}\label{sec:U}
Recall that the module $V$ is generated over $R$ by a set $\Cc$ and that the module $M$ is generated over $A$ as a left-right module by $\wh \Cc\cup \Xx$, where $\wh \Cc=\{\wh c\colon c\in \Cc\}$ and $\Xx=\{x_1,\dots,x_m\}$. 
Let $\wh V^\ast$ denote the $K$-vector space dual to $\wh V$. Then $\wh V^{\ast}$ is the {\em direct product}  of the 1-dimensional $K$-vector spaces  which are generated by elements of  $\wh \Cc^{\ast}\cup \Xx^{\ast}$, where
 $\wh \Cc^{\ast}=\{\wh c^{\,\ast}\colon c\in \Cc\}$ and  $\Xx^{\ast}=\{x_j^{\ast}\colon x_j\in \Xx\}$.
Furthermore, we have 
\[
\wh V^\ast=\bigoplus_{1\le i,j\le m}
\wh V^{\ast\,(ij)},\quad\text{where}\quad
\wh V^{\ast\,(ij)}:=\bigl(\wh V^{(ij)}\bigr)^{\ast}. 
\]
In $\wh V^{\ast}$, the idempotent  $1_j$, $j=1,\dots, m$, acts from the left as the identity on $\wh V^{\ast\,(ji)}$ and as $0$ on $\wh V^{\ast\,(li)}$, $l\ne j$, and acts from the right as the  identity on $\wh V^{\ast\,(ij)}$ and as $0$ on $\wh V^{\ast\,(il)}$, $l\ne j$. 
Fix an integer $n \ge 2$ (it will correspond to half the dimension of the Liouville domain $X_0$
in Section~\ref{sec:leg-algebra}).
We endow $\wh V^{\ast}$ with a grading defined as follows:
\begin{align}
&|{\wh c}^{\,\ast}|=n-3-|\wh c\,|, \text{ for }\wh c\in\wh \Cc,\\ 
&|x_j^{\ast}|=n-3,\text{ for }j=1,\dots,m.
\end{align}

Define  
\begin{equation}\label{eq:defM*}
M^{\ast}:=A\mathop{\otimes}\limits_R\wh V^{\ast}\mathop{\otimes}\limits_R A.
\end{equation}
We define the tensor algebra
\begin{equation}\label{eq:U(0,0)}
U(0,0):= {\bigoplus}_{k\geq 0}\, (M\oplus M^\ast)^{{\motimes\limits_A}^k}.
\end{equation}
as a certain completion of the free associative algebra generated by $\Cc\cup \wh\Cc\cup \wh\Cc^\ast\cup\XX\cup\XX^\ast$.  Its elements are infinite series in the $(\wh\Cc^{\ast}\cup\Xx^{\ast})$-variables with polynomial coefficients in the $(\Cc\cup\wh\Cc\cup\Xx)$-variables (compare to the above definitions of $M$ as an infinite sum and $M^{\ast}$ as an infinite product). For simplicity of notation and since all completions we consider below are of the same type as that of \eqref{eq:U(0,0)}, we will, as in that equation, suppress completions from notation and use the direct sum symbol in decompositions.
Consider the extension $U$ of $U(0,0)$ obtained by adding multiplicative generators
$\hbar,\hbar^{-1}$ and $\sigma,\sigma^{-1}$ of degrees 
\[
|\sigma|=1,\quad |\sigma^{-1}|=-1,\quad |\hbar|=n-3,\quad|\hbar^{-1}|=-(n-3)
\] 
which satisfy the following relations: 
\begin{equation*}
\sigma\sigma^{-1}=\sigma^{-1}\sigma=1,\quad \hbar\hbar^{-1}=\hbar^{-1}\hbar=1,\quad
\sigma\hbar=(-1)^{n-3}\hbar\sigma,
\end{equation*}
and
\begin{align*}
\sigma u&=\begin{cases}
(-1)^{|u|}u\sigma,& u\in\Cc;\cr
(-1)^{|u|+1}u\sigma,&  u\in\wh \Cc\cup\XX\cup\wh \Cc^\ast\cup \XX^\ast,
\end{cases}\\
\hbar u&=\begin{cases}
(-1)^{|u|(n-3)}u\hbar,& u\in\Cc\cup\Xx\cup\wh\Cc;\cr
(-1)^{(|u|+n-3)(n-3)}u\hbar,&  u\in\wh \Cc^\ast\cup\XX^\ast. 
\end{cases}
\end{align*}

The algebra $U$ can be decomposed according to tensor type: $U=\bigoplus_{p,q\geq 0}U^p_q$,
where $U^p_q\subset U$ is the $K$-subspace generated by monomials with $q$ factors from $M$ and $p$ from $M^\ast$. 
We decompose further:
$U^p_q=\bigoplus_{s,h\in\Z} U^p_q(s,h)$, where $U^p_q(s,h)\subset U^{p}_q$ is the subspace spanned by all monomials of total $\sigma$-degree $s$ and total $\hbar$-degree $h$.
In particular, $U^0_0(0,0)=A$, $U^0_1(0,0)=M$, and  $U^1_0(0,0)=M^\ast$. We also write $U^{p}:=\bigoplus_{q\ge 0} U^{p}_{q}$, $U^{+}:=\bigoplus_{p>0} U^{p}$, $U_{q}:=\bigoplus_{p\ge 0}U^{p}_{q}$, $U_q(s,h):=\bigoplus_{p\ge 0} U^{p}_{q}(s,h)$, and $U^{p}(s,h):=\bigoplus_{q\ge 0}U^{p}_{q}(s,h)$.

For $X\in U^p_q(s,h)$, let 
\[
\st(X):=p+h\quad\text{and}\quad \ju(X):=p+q+s. 
\]
In analogy with the definition of $M^\diag$ in Section \ref{sec:MMcyc}, we define $U^\diag$. Further, define $U^\cyc$ as the quotient of $U^{\diag}$ obtained by dividing by the following {\em graded cyclic permutation} rule: if $a_1, \dots,a_k\in\Cc\cup\wh\Cc\cup\Xx\cup\wh\Cc^{\ast}\cup\Xx^{\ast}$ and if
$a_1\dots a_k\in U^{\diag}$ then $a_2\dots a_ka_1\in U^{\diag}$, and we identify
$a_1\dots a_k\in U^{\diag}$ with $\eps_1\eps_2\eps_3\;a_2\dots a_ka_1\in U^{\diag}$, where
\begin{align*}
&\eps_1=(-1)^{|a_1||a_2\dots a_k|},\\
&\eps_2={\begin{cases}
1&\text{if } a_1\in\Cc\cup\{\hbar,\hbar^{-1}\},\cr
(-1)^{\ju(a_2\dots a_k)}&\text{if }a_1\in\wh\Cc\cup\Xx\cup\wh\Cc^{\ast}\cup\Xx^{\ast}\cup\{\sigma,\sigma^{-1}\}, \cr
\end{cases}}\\
&\eps_3={\begin{cases}
1 &\text{if }a_1\in\Cc\cup\wh\Cc\cup \Xx\cup\{\sigma,\sigma^{-1}\},\cr
(-1)^{\st(a_2\dots a_k)}&\text{if }a_1\in\wh\Cc^{\ast}\cup\Xx^{\ast}\cup\{\hbar,\hbar^{-1}\}.
\end{cases}}
\end{align*}
We will often think of cyclic monomials as written on a circle $S^{1}$ with the first letter at $1\in S^{1}$. This allows us to speak about generators in a cyclic monomial as ordered in the (counter) clockwise direction staring from some fixed generator in the monomial. In this language, the above permutations correspond to rotations of the circle.

The vector space $U^{\cyc}$ inherits the decomposition by tensor type:  
\[
U^{\cyc}=\bigoplus_{p,q\geq 0} (U^{p}_q)^{\cyc}=\bigoplus_{p,q\geq 0}\bigoplus_{s,h\in\Z}\,\lb U^{p}_q(s,h)\rb^{\cyc}.
\]
If the monomials $X$ and $X'$ differ by a graded cyclic permutation, i.e.~if $\llb X\rrb=
\llb X' \rrb$, then, for $z\in\{\sigma,\sigma^{-1},\hbar,\hbar^{-1}\}$, $\llb Xz\rrb =\llb X'z\rrb$  and  
$\llb zX\rrb =\llb zX'\rrb$. We thus define $z\llb X\rrb:=\llb zX\rrb$ and $\llb X\rrb z:=\llb Xz\rrb$, and note that this definition does not depend on choice of representative $X$.
  
A monomial $X\in\lb U^{p}_q(s,h)\rb^{\cyc}$ is {\em balanced} if 
$h=-1$ and $p-q-s=0$.
The subspace of $U^\cyc$ generated by balanced monomials will be denoted by $U^\bal$. It contains balanced versions $M^\bal$ and $\lb M^{\ast}\rb^{\bal}$ of $M^\cyc$ and $\lb M^{\ast}\rb^{\cyc}$, respectively. More precisely, 
\begin{align*}
M^\bal&:=\lb U^0_1(-1,-1)\rb^\cyc= U^0_1\cap U^\bal,\\ 
\lb M^{\ast}\rb^{\bal}&:=\lb U^1_0(1,-1)\rb^\cyc= U^1_0\cap U^\bal.
\end{align*}

\subsubsection{The operator $\Ss\colon U^{\cyc}\to U^{\cyc}$}
Consider a monomial $X\in \lb U^p_q(s,h)\rb^\cyc$ and fix a factor $c\in\Cc$ in $X$. Then $X=\llb X_c'cX_c''\rrb$ for monomials $X_c',X_c''\in U$. Define    
\[
\Ss(X)=\sum\limits_{c}\llb X_c'\sigma^{-1}\wh cX_c''\rrb\in \lb U^p_{q+1}(s-1,h)\rb^\cyc,
\]
where the sum is taken over all factors $c\in \Cc$ in $X$. It is straightforward to check that this gives a well defined operation and by definition $\Ss(U^{\bal})\subset U^{\bal}$.
In particular, if $X=\llb c_1\dots c_m \wh a^\ast\rrb\in \lb U^1_0(0,0)\rb^\cyc$ then $\Ss(X)\in\lb  U^0_{1}(-1,0)\rb^\cyc $ satisfies
\begin{align*}
\Ss(c_1\dots c_m\wh a^{\ast})&=\sigma^{-1}\sum_{j=1}^{m} (-1)^{|c_1\dots c_{j-1}|}c_1\dots c_{j-1}\wh c_j c_{j+1}\dots c_m\wh a^\ast\\
&=\sigma^{-1} S(c_1\dots c_m)\wh a^\ast, 
\end{align*}
where $S$ is the operator $A\to M$ defined in \eqref{eq:Sop} above. Furthermore, 
 \[
\Ss^k(c_1\dots c_m \wh a^{\ast})=k!\sigma^{-k}\sum_{j_1<\dots<j_k}^{m} (-1)^{\sum\limits_{l=1}^k|c_1\dots c_{j_l-1}|}  c_1\dots\wh c_{j_1}\dots \wh c_{j_k}\dots c_m \wh a^{\ast}.
\]

\subsubsection{Symmetry breaking and contraction operations}\label{sec:convol}
We next effectively break the cyclic symmetry of $U^\bal$. An {\em excited} monomial is a cyclic monomial $X\in U^{\bal}$ with one of its factors from $\wh\Cc\cup\wh\CC^\ast\cup \Xx\cup\Xx^\ast$ distinguished. We call the distinguished factor of an excited monomial $X$ the \emph{excited} generator of $X$. The $K$-vector space generated by excited monomials from $U^\bal$ will be denoted by $\bU^\bal$. 
Note that $\bU^{\bal}$ inherits decompositions from $U^{\bal}$, $\bU^{\bal}=\bigoplus_{p,q} (\bU^{p}_{q})^{\bal}$, etc. We set $(\bU^0_0)^{\bal} = 0$. Define the {\em excitation operator} $\Ee:U^\bal\to \bU^\bal$ to be the linear map which maps a monomial $X$ the sum of all its excitations.

Note that $\bU^{\bal}$ contain excited versions $\bM^\bal=\lb\bU^{0}_{1}\rb^{\bal}$ and $\lb\bM^{\ast}\rb^{\bal}=\lb\bU^{1}_{0}\rb^{\bal}$ of $M^\bal$ and $\lb M^{\ast}\rb^{\bal}$, respectively, which are isomorphic to their non-excited counterparts. (In the context of $\bU^{\bal}$ we think of $M^{\bal}$ and $\lb M^{\ast}\rb^{\bal}$ as generated by monomials in $\lb U^0_1\rb^{\bal}$ and $\lb U^{1}_0\rb^{\bal}$ {\em without} excited generator.)

We will use the following notation to identify excited generators in monomials: we underline the excited generator in a monomial in $U^{\bal}$ and write e.g.~$\ux$, $\ux^\ast$ $\wh\uc$, and $\wh\uc^\ast$.

We next define a product on $\bU^\bal$. Consider two excited cyclic monomials $X,Y\in \bU^\bal$. Fix a factor $u\in \wh\Cc^\ast\cup\XX^\ast$ from $X$ and a factor $v\in \wh \Cc \cup\XX $ from $Y$. Write $X=\eps_1\llb X_u'u\rrb$ and $Y=\eps_2 \llb v Y_v'\rrb$, where $\eps_1$ is the sign which arises from cyclically rotating $X$ so that $u$ is its last factor and where $\eps_2$ is the sign which arises from cyclically rotating $Y$ so that $v$ is its first factor. Define the {\em contraction} operation $\mstar\colon\bU^{\bal}\otimes\bU^{\bal}\to\bU^{\bal}$,
\begin{equation}\label{eq:contraction}
X\star Y:=\sum\limits_{(u,v)}\eps_1\eps_2 u(v)\llb X_u'Y_v'\rrb,
\end{equation}
where  the sum is taken over all pairs $(u,v)$ of factors $u\in \wh\Cc^\ast\cup\XX^\ast$ from $X$ and $v\in \wh \Cc \cup\XX $ from $Y$, and where
\[
u(v)=
\begin{cases}
\hbar&\text{if exactly one of $u$ and $v$ is excited, $u=v^\ast$, and $v\notin \Xx$ ,}\\
(-1)^{n-1}\hbar&\text{if $u=\ux_j^\ast$ and $v=x_j$, for some $j=1,\dots,m$,}\\
\hbar &\text{if $u=x_j^{\ast}$ and $v=\ux_j$, for some $j=1,\dots,m$,}\\
0&\text{otherwise}.
\end{cases}
\]

If $Y\in \lb \bU^1_q\rb^\bal$ and $X\in\bU^{\bal}$ then we define, for $j=1,\dots, q$, the {\em partial contraction} $X\mstar_j Y$ as follows:  
\[
X\mstar\limits_j Y:=\sum\limits_{u}\eps_1\eps_2u(v_j)\llb X_u'Y_{v_j}'\rrb,
\]
where the sum is taken over all $u\in \wh\Cc^\ast\cup\XX^\ast$ from $X$ and where $v_j$ is the $j^{\rm th}$ factor from 
$\wh \Cc \cup\XX $ in $Y$ counting {\em counter-clockwise} from the unique factor from $\wh\Cc^\ast\cup\XX^\ast$.  Similarly, if  $X\in \lb \bU^p_1\rb^\bal$ and $Y\in\bU^{\bal}$ then we define, for $i=1,\dots,p$ the {\em partial contraction} $X\mstar^i Y$:  
\[
X\mstar\limits^i Y:=\sum\limits_{v}\eps_1\eps_2u_i(v)\llb X_{u_i}'Y_v'\rrb,
\]
where the sum is taken over all $v\in \wh\Cc \cup\XX $ from $Y$ and where $u_i$ is the $i$-th factor from  $\wh \Cc^\ast \cup\XX^\ast $ in $X$ counting {\em clockwise} from the unique factor from $\wh\Cc \cup\XX $. 
\begin{rmk}
Let $X,Y\in \bU^\bal$ then $X\star Y\in \bU^\bal$. Further if $X\in\lb\bU^{1}_{q}\rb^{\bal}$ then $X\mstar_j Y\in\bU^{\bal}$, and if $Y\in\lb\bU^{p}_1\rb^{\bal}$ then $X\mstar^{i}Y\in\bU^{\bal}$.  
\end{rmk}	
Define the {\em bracket} or {\em commutator} $[\,,]\colon \bU^{\bal}\otimes\bU^{\bal}\to\bU^{\bal}$ as 
\begin{equation}\label{eq:bracket}
[X,Y]=X\star Y -(-1)^{|X|  |Y|}\,Y\star X.
\end{equation}
\begin{lma}\label{l:Jacobi}
The bracket satisfies the graded Jacobi identity:
\[
(-1)^{|X||Z|}[X,[Y,Z]]+(-1)^{|Z||Y|}[Z,[X,Y]]+(-1)^{|Y||X|}[Y,[Z,X]]=0.
\]
\end{lma}

Let 
\[
\tau\bM^{\bal}=\bM^{\bal}\oplus M^{\bal}\quad\text{and}\quad 
\lb\tau\bM^{\ast}\rb^{\bal}=\lb\bM^\ast\rb^\bal\oplus \lb M^\ast\rb^\bal.
\]
We associate with $X\in \lb\bU^p_1\rb^\bal$ a multi-linear operation $X^\downarrow\colon\lb \tau\bM^{\bal}\rb^{\otimes^p}\to\tau\bM^{\bal}$ defined by
\begin{equation}\label{eq:down}
X^\downarrow (A_1\otimes\dots\otimes A_p)= \lb\dots\lb\lb X\mstar^1 A_1\rb\mstar^1 A_2\rb\mstar^1 \dots\rb\mstar^1 A_p.
\end{equation}
Note that if $A_1,\dots, A_p\in \bM^\bal$ then $X^\downarrow (A_1\otimes\dots\otimes A_p)\in \bM^\bal$, that if exactly one of the arguments lies in $M^\bal$ then the image lies in $M^\bal$, and that the operation vanishes on $p$-tuples with more than component in $M^\bal$.

Similarly, $Y\in \lb\bU_q^1\rb^\bal$ defines a multi-linear operation 
$Y^\uparrow\colon\lb\lb\tau\bM^{\ast}\rb^\bal\rb^{\otimes^q}\to\lb\tau\bM^\ast\rb^\bal$,
\begin{equation}\label{eq:up}
Y^{\uparrow}(B_1\otimes\dots\otimes B_q)= B_q\mstar_1\lb\dots\mstar_{1} \lb B_{2}\mstar_{1}\lb B_1\mstar_1 Y\rb\rb\dots\rb,
\end{equation} 
if $B_1,\dots, B_q\in \lb\bM^\ast\rb^\bal$ then $Y^\uparrow (B_1\otimes\dots\otimes B_q)\in \lb\bM^\ast\rb^\bal$, if exactly one  of the arguments lies in $\lb M^\ast\rb^\bal$ then the image lies in $\lb M^\ast\rb^\bal$, and the operation vanishes on $q$-tuples with more than one component in  $\lb M^\ast\rb^\bal$.

In the special case $q=2$, we associate with  $Y\in \lb \bU^1_2\rb^\bal$  two more operations, 
\begin{equation}\label{eq:updown}
Y^{\uda} \colon \lb\tau\bM^\ast\rb^\bal\otimes \tau \bM^\bal \to \tau\bM^\bal
\end{equation}
and 
\begin{equation}\label{eq:downup}
Y^{\dua} \colon  \tau\bM^\bal \otimes\lb \tau\bM^\ast\rb^\bal\oplus \to\tau\bM^{ \bal},
\end{equation}
where
\begin{align*}
Y^{\uda}(X,Z)&= \lb X\mstar_1 Y\rb\star Z,\quad X\in \lb\tau\bM^\ast\rb^\bal,\; Z\in 
\tau\bM^\bal,\\ 
Y^{\dua}(Z,X)&= \lb X\mstar_2Y\rb\star Z ,\quad Z\in\tau\bM^\bal,\; X\in \lb\tau\bM^\ast\rb^\bal.
\end{align*}
Similarly, $X\in\lb\bU^2_0\rb^\bal$  induces an operation $X^{\circlearrowleft}\colon\tau \bM^\bal \to \lb \tau\bM^
\ast\rb^\bal$,
\begin{equation}\label{eq:circ}
X^{\carrl}(Y)=  X\star Y.
\end{equation}
We observe that $X^{\carrl}(\bM^\bal)\subset\lb\bM^\ast\rb^\bal$ and   $X^{\carrl}(M^\bal)\subset \lb M^\ast\rb^\bal$.

\subsubsection{The Hamiltonian}\label{sec:Hamiltonian}

We define special elements $H^1_q\in\bU^\bal$ which derive from the differential on $A$ and which, in the spirit of SFT, see \cite{SFT}, are called {\em Hamiltonians}.
\begin{dfn}\label{dfn:hamiltonian1}
\begin{align*}
&(h^1_1)' :=\sum\limits_{c\in\Cc}\llb S(dc)\hbar^{-1} \wh c^\ast\rrb \\&
(h^1_1)'' := 
\sum\limits_{i,j=1}^m\sum\limits_{c\in\Cc_{ij}}\lb \llb c x_j\hbar^{-1}\wh c^\ast\rrb-\llb x_i c \hbar^{-1}\wh c^\ast\rrb\rb 
 \\&
 h^1_1 := (h^1_1)' + (h^1_1)'' \\&
 h^{1}_2:=\frac{1}{2!} \Ss (h^{1}_1)' +  \Ss (h^{1}_1)'' +  (-1)^{n-1}\sigma^{-1}\sum\limits_{j=1}^m\llb x_jx_j \hbar^{-1}x_j^{\ast} \rrb
\\
&h^1_q:=\frac1{q!}\Ss^{q-1} h^1_1,\quad \text{ for }q>2,\\
&H^1_q:=\Ee(h^1_q),\;\; \cH^1=\sum\limits_{q\geq 1}\cH^1_q.
\end{align*}
\end{dfn}
 
We use the Hamiltonian $H^1_1$ to define a differential on $\tau\bU^\bal:=\bU^{\bal}\oplus M^{\bal}\oplus (M^{\ast})^{\bal}$. First define the {\em exterior differential} $d\colon \tau\bU^\bal\to \tau\bU^\bal$ as follows. Consider the map $d\colon U\to U$ which acts as the algebra differential $d\colon A\to A$ on generators of $U$ from $A$, which maps all other generators to $0$, and which satisfies the graded Leibniz rule. Then $d$ preserves $U^{\diag}$, descends to $U^{\cyc}$ where it preserves $U^{\bal}$. It then induces a map $d\colon \tau\bU^{\bal}\to\tau\bU^{\bal}$ as follows. The map is already determined on $M^{\bal}\oplus (M^{\ast})^{\bal}$ and
if $\llb\ul u X\rrb$ is a monomial in $\bU^{\bal}$, where $u\in\Cc\cup\Xx\cup\wh\Cc\cup\wh\Cc^{\ast}\cup\Xx^{\ast}$ then $d(\llb\ul u X\rrb)=(-1)^{|u|}\llb \ul u(dX)\rrb$. 
Then $d^2=0$ and we observe that $H^1$ satisfies the following {\em master equation}
\begin{equation}\label{eq:master1}
d\cH^1+\cH^1\star \cH^1=0.
\end{equation}
Since $|H^{1}|=-1$, $\cH^1\star \cH^1=\frac12[\cH^1, \cH^1]$ and \eqref{eq:master1} can be rewritten in Maurer-Cartan form:
\begin{equation}
\label{eq:master2}
d\cH^1+\frac12[\cH^1, \cH^1]=0.
\end{equation}

Decomposing $\cH^1=\sum_{q\geq 1}\cH^1_q$, \eqref{eq:master1} is equivalent to the following sequence of identities  
for the $\lb \bU^1_q\rb^\bal$-components of $\cH^{1}\star \cH^{1}$:
\begin{equation}\label{eq:I_q}
d\cH_q^{1}+\sum_{k=1}^q \sum_{j=1}^k \cH^1_{q-k+1}\mstar_j \cH^1_{k}=0,\quad q=1,2,3,\dots
\end{equation} 
The first identity in this sequence is
\begin{equation}
\label{eq:master0}
d\cH^1_1+\cH^1_1\star \cH^1_1=0,
\end{equation}
\begin{dfn} Define the differential $d_U\colon\tau\bU^{\bal}\to\tau\bU^{\bal}$,
\[
d_U X:=dX+[H^1_1,X].
\]
(Note that \eqref{eq:master0} is equivalent to $d_U^{\,2}=0$.)
\end{dfn}
The differential interacts with the bracket in the following way.
\begin{lma}\label{lm:star}
For $X,Y\in\tau\bU^{\bal}$,
\[
d_U([X,Y])=[d_UX, Y]+(-1)^{|X|}[X, d_UY].
\] 
In particular, the operation $[\,,]$ descends to the homology $H(\tau\bU^\bal,d_U)$.
\end{lma}

We observe that  $d_U$ leaves $M^\bal$, $\bM^\bal$, $\lb M^\ast\rb^\bal$, and $\lb\bM^\ast\rb^\bal$ invariant. For simple notation we write $d_{\tau M}$  for the restriction of $d_U$ to any one of these subspaces. Consider the isomorphisms  $\beta\colon M^\cyc\to M^\bal$ and $\ul\beta\colon M^{\cyc}\to\bM^{\bal}$, 
\begin{equation}\label{eq:beta}
\beta(X)=X\hbar^{-1}\sigma^{-1}\;\text{ and }\;{\ul\beta}=\Ee\circ\beta,\;\text{ respectively}.
\end{equation}
\begin{lma}
The maps
 $\beta\colon(M^\cyc, d_M)\to (M^\bal, {d_\tau M})$  and  ${\underline\beta}\colon(M^\cyc, d_M)\to (\bM^\bal, d_{\tau M})$
are chain isomorphisms.
\end{lma} 
The following result is a consequence of Lemma \ref{lm:star}.
\begin{cor}
If $X\in \bU^{\bal}$ satisfies $d_{U } X=0$ and has the correct tensorial type to define one of the operations \eqref{eq:down} -- \eqref{eq:circ}, see Section \ref{sec:convol}, then this operation is a chain map. If furthermore $X=d_U Y$ for some $Y \in \bU^\bal$ then the operation is chain homotopic to $0$.
\end{cor}

\subsection{Products}\label{sec:operations}
We introduce product operations on $\bM^\bal$ and $\lb\bM^\ast\rb^\bal$.
\subsubsection{The product on $H(\lb\bM^\ast\rb^\bal)$}
Recall the operation  $\lb \cH^1_2\rb^{\uparrow}\colon  
  \lb\tau\bM^\ast\rb^\bal \otimes  \lb \tau\bM^\ast\rb^\bal \to
  \lb\tau\bM^\ast\rb^\bal$ defined by the formula
\begin{equation*} 
 \lb \cH^1_2\rb^{\uparrow}(X,Y)=Y\star\lb X\mstar\limits_1\cH^1_2\rb,
\end{equation*}
see \eqref{eq:up}.
Define $X\diamond Y:=(-1)^{|X|} \lb \cH^1_2\rb^{\uparrow}(X,Y)$. Then $\lb\bM^\ast\rb^\bal\diamond \lb\bM^\ast\rb^\bal\subset \lb \bM^\ast\rb^\bal$.
 \begin{prp}\label{prop:dual-product}
The operation
$\diamond\colon\lb\bM^\ast\rb^\bal \otimes\lb\bM^\ast\rb^\bal\to\lb\bM^\ast\rb^\bal$ descends to an associative product on the homology
$H(\lb\bM^\ast\rb^\bal,d_{\tau M})$ which satisfies the following  graded commutativity relation for the grading shifted by $1$:
\[
X\diamond Y=(-1)^{(|X|+1)(|Y|+1)}Y\diamond X.
\]
The  homology class of the cycle $E:=(-1)^{n-1}\sum_{i=1}^k \hbar^{-1}\ux_{\,i}^\ast  \sigma   \in  
(\bM^{\ast})^\bal$  is the unit for this product  
on $H(\lb\bM^\ast\rb^\bal)$.  In fact, $E$ is a unit for $\diamond$ already on the chain level, before passing to homology.
\end{prp}
Consider the second and third identities satisfied by the Hamiltonian $\cH^{1}$:
\begin{align}\label{eq:I_2}
&d\cH^{1}_{2}+[\cH^1_1,\cH^1_2]=0,\\\label{eq:I_3}
&d\cH^{1}_{3}+[\cH^1_1,\cH^1_3]+ \cH^1_2\mstar _1 \cH^1_2 + \cH^1_2\mstar _2 \cH^1_2=0.
\end{align}
Here \eqref{eq:I_2} is the chain map equation for $\lb \cH^1_2\rb^{\uparrow}\colon\lb \bM^\ast\rb^\bal\otimes \lb \bM^\ast\rb^\bal\to\lb \bM^\ast\rb^\bal$, and \eqref{eq:I_3} shows that this product is associative on the homology level.
Indeed,
\begin{align*}
(X\diamond Y)\diamond Z
&= (-1)^{|Y|+1}Z\mstar ((Y\star (X\mstar_1H^1_2))\mstar_1 H_1^2)\\
&=(-1)^{|Y|}Z\mstar (Y\mstar_1 (X\mstar_1 (H^1_2\mstar_2 H_1^2)))\\ 
&=(-1)^{|Y|+|X|}(Z\mstar(Y\mstar_1 H^1_2))\mstar (X\mstar_1  H_1^2)\\
&= X\diamond (Y \diamond Z).
\end{align*}
The unit property can be checked directly if one  observes that the sum $G$ of all terms in $\cH^1_2$ with excited generator from $\wh\Cc^{\ast}\cup\Xx^{\ast}$ and which contain at least one $x_j$-factor, $j=1,\dots,m$ is 
\[
G= \sigma^{-1}\lb(-1)^{n-1}\sum_j\llb x_j x_j\hbar^{-1} \ux^\ast_j\rrb +\sum_{i,j}\sum\limits_{c\in\Cc_{ij}}
 \llb x_i\wh c\hbar^{-1}\wh \uc^{\,\ast}\,\rrb +\llb\wh c x_j\hbar^{-1}\wh \uc^{\, \ast} \rrb\rb.
\]
Hence,  taking into account that $|E|=1$  we have for any $Y\in \lb \bM^\ast\rb^\bal$  
\begin{align*}
&E\diamond Y=Y\star \left( \sum_{i=1}^m (-1)^{n-1}\hbar^{-1}\ux_i^\ast \sigma  \mstar_1 G\right)\\&
%=Y\star \sum_{i=1}^k \hbar^{-1}\ux_i^\ast  \mstar^1 \lb\sum_j\llb x_j x_j\hbar^{-1} \ux^\ast_j\rrb +(-1)^{n-1}\sum_{i,j}\sum\limits_{c\in\Cc_{ij}}
 %\llb x_i\wh c\hbar^{-1}\wh \uc^{\,\ast}\,\rrb +\llb\wh c x_j\hbar^{-1}\wh \uc^{\, \ast} \rrb\rb
&\\&=
Y\star\lb  \sum_j(-1)^{n-1}\llb x_j\hbar^{-1}\ux^\ast_j\rrb + \sum\limits_{c\in\Cc }
 \llb \wh c\hbar^{-1}\wh \uc^{\,\ast}\,\rrb \rb=Y
\end{align*}
 
To verify the commutativity relation, let $X$ be a $d_M$-cycle in $\bM^{\bal}$. Applying the $\Ss$-operator to the equation $dX+ [H^1_1, X ]=0$ gives $d_U(\Ss X)=X\mstar H^{1}_2$ and thus, on the level of homology, $X\mstar_1 H^1_2+X\mstar_2 H^1_2=0$ and we calculate 
\begin{align*}
&X\diamond Y=(-1)^{|X|}Y\mstar (X\mstar_1 H^1_2)=(-1)^{|X|+1}Y\mstar(X\mstar_2 H^1_2)\\
&=
(-1)^{|X|+1+|X||Y|}X\mstar(Y\mstar_1H^1_2)=
(-1)^{(|X|+1)(|Y|+1)}Y\diamond X.
\end{align*}
\subsubsection{The product on $H(\bM^\bal)$}\label{sec:prod}
If there exists an element 
\[
\cH'=\sum_{p\geq 2} \cH^p\;\in\;\bigoplus\limits_{p\geq 2} \lb \bU^p\rb^\bal
\]
such that $\cH= \cH^1+\cH'=\sum_{p\geq 1} \cH^p$  satisfies the master equation
\begin{equation}\label{eq:fullmaster}
d\cH+\cH\star \cH=d\cH+\frac12[\cH,\cH]=0,
\end{equation} 
then we say that $\cH'$ is {\em compatible} with $\cH^1$ and we call $\cH=\cH^1+H'$ the {\em full (rational) Hamiltonian}. As we shall see later, when $A$ is a Legendrian homology algebra then we can always find $\cH'$ compatible with $\cH^1$ using  a relative SFT formalism associated with parallel copies of the Legendrian submanifold in the spirit of \cite{E:rsft}.

Assume now that $\cH$, and in particular $\cH^2_0$, which satisfies \eqref{eq:fullmaster} exists and is given. Then we define the operation $\boxdot\colon\tau\bM^{\bal}\otimes\tau\bM^{\bal}\to\tau\bM^{\bal}$, 
\[
X\boxdot Y= \lb\cH^2_0\mstar_1 \cH^1_2\rb^{\downarrow}(X,Y)=((\cH^2_0\mstar_1 \cH^1_2)\mstar^1 X)\mstar Y.
\] 
Note that if $X,Y\in\bM^{\bal}\subset\tau\bM^{\bal}$ then $X\boxdot Y\in\bM^{\bal}$.
\begin{thm}\label{thm:product-properties}
${ }$
\begin{enumerate}
\item The operation $\boxdot\colon\bM^{\bal}\otimes\bM^{\bal}\to\bM^{\bal}$ descends to an operation on the homology $H(\bM^\bal,d_{\tau M})$ where it can be expressed by the following three equivalent expressions: 
\begin{align*}
X\boxdot Y=&-\lb\cH^2_0\mstar_2 \cH^1_2\rb^{\downarrow}(X,Y) 
 =(-1)^{|X|}\lb H^1_2\rb^{\uda}\lb \lb H^2_0\rb^\carrl(X),Y\rb\\
 =&\;(-1)^{|X||Y|+|X|+1}\lb H^1_2\rb^{\dua}\lb X, \lb H^2_0\rb^{\carrl}(Y)\rb.
\end{align*}
Equivalently, the following diagram commutes on the homology level,
   \[
  \xymatrix{&\lb \bM^\ast\rb^\bal\otimes \bM^\bal\ar[dr]^{(\cH^1_2)^{\wt\dua}}&\\
  \bM^\bal\otimes \bM^\bal\ar[ur]^
  {(\cH^2_0)^\carrl\otimes\Id\;}\ar[rr]^{\boxdot\;\simeq\; -(\cH^2_0\mstar_2 \cH^1_2)^{\downarrow}}\ar[dr]_{\Id\otimes(\cH^2_0)^\carrl\;\;}&& \bM^\bal\\
  & \bM^\bal \otimes\lb \bM^\ast\rb^\bal\ar[ur]_{(\cH^1_2)^{\wt\uda}}&}
  \]
  where 
\begin{align*}
(\cH^1_2)^{\wt\uda}(X,Y)&:=(-1)^{|X|}(\cH^1_2)^{ \uda}(X,Y),\\ 
(\cH^1_2)^{\wt\dua}(Y,X)&:=(-1)^{|Y||X|+|Y|+1}\linebreak(\cH^1_2)^{\dua}(Y,X),
\end{align*}
 for $X\in \lb  \bM^\ast\rb^\bal$ and $Y\in \bM^\bal$.
 \item The operation $\boxdot$ is  associative and graded commutative on homology, i.e.~if $X,Y\in H(\bM^{\bal})$ then    
$X\boxdot Y=(-1)^{|X||Y|} Y\boxdot X.$  
\end{enumerate}
\end{thm}

The above result is a consequence of the master equation \eqref{eq:fullmaster} which in particular implies the identities 
\begin{align}\label{eq:H20}
&d\cH^{2}_{0}+ [\cH^1_1,\cH^{2}_0]=d\cH^{2}_{0}+\cH^2_0\star \cH_1^1=0,\\\label{eq:H21}
&dH^{2}_{1}+[H^1_1,H^2_1]+H^2_0\mstar_1 H^1_2+ H^2_0\mstar_2H^1_2=0.
\end{align}

A straightforward calculation shows that $d_U(H^{2}_{0}\mstar_1 H^{1}_2)=(d_U H^{2}_0)\mstar_1 H^{1}_2- H^{2}_0\mstar_1(d_U H^{1}_{2})$. 
By \eqref{eq:H20}, $d_{U}\cH^2_0=0$, by \eqref{eq:I_2}, $d_{U}\cH^{1}_{2}=0$, and therefore $\cH^2_0\mstar_1 \cH^1_2$ is a cycle in $\lb \bU^2_1\rb^\bal$. Similarly, $-\cH^2_0\mstar_2 \cH^1_2$ is a cycle and \eqref{eq:H21} implies that their homology classes coincide. This shows that the product descends to homology and implies the first equation of part 1 of the theorem. The alternative expressions for $\boxdot$ in part 1 are immediate from the definitions of the operations $ \lb H^2_0\rb^\carrl$, $\lb H^1_2\rb^{\dua}$, and $\lb H^1_2\rb^{\uda}$.

\begin{figure} 
\labellist
\small\hair 2pt
\pinlabel	$(\mathrm{I})$ at -30 460
\pinlabel $+$ at 90 460
\pinlabel $+$ at 208 460
\pinlabel $+$ at 316 460
\pinlabel $+$ at 428 460
\pinlabel $=\;0$ at 550 460
\pinlabel	$(\mathrm{II})$ at -30 300
\pinlabel $+$ at 72 300
\pinlabel $+$ at 175 300
\pinlabel $+$ at 264 300
\pinlabel $+$ at 360 300
\pinlabel $+$ at 450 300
\pinlabel $=\;0$ at 550 300
\pinlabel $\stackrel{(\mathrm{I})}{\simeq}$ at 170 120
\pinlabel $\stackrel{(\mathrm{II})}{\simeq}$ at 400 120
\endlabellist
\centering
\includegraphics[width=.5\linewidth]{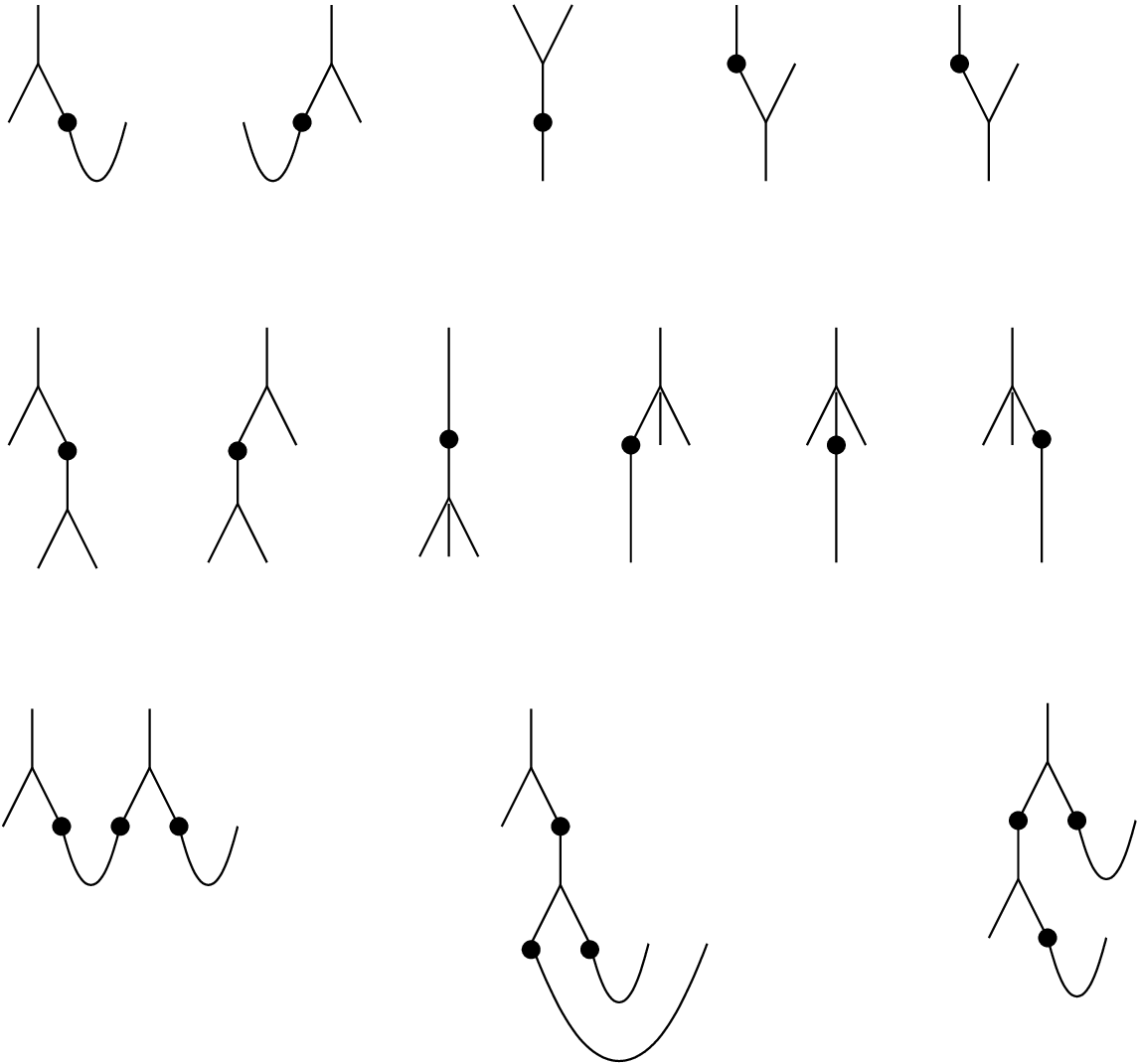}
\caption{A graphical proof of associativity. A tree with $p$ upwards and  $q$ downwards  $1$-valent vertices represents $\cH^{p}_q$. The dot represents the $\mstar$-operator. Equations $(\mathrm{I})$ and $(\mathrm{II})$ are consequences of $d\cH+\cH\mstar \cH=0$ (terms from $d\cH$ are not shown in the figures).}
\label{fig:calc}
\end{figure}
To verify the  commutativity we observe that given a cycle $X\in \bM^\bal$,
using the homological identity $(H^2_0\mstar X)\mstar_1 H^1_2=-(H^2_0\mstar X)\mstar_2 H^1_2$, we calculate
\begin{align*} 
&X\boxdot Y =((\cH^2_0\mstar_1 \cH^1_2)\mstar^1 X)\star Y=
(-1)^{|X|} ((H^2_0\mstar X)\mstar_1 H^1_2)\star Y\\
&=(-1)^{|X|+1} ((H^2_0\star X)\mstar_2H^1_2)\mstar Y=
(-1)^{|X||Y|+1 }((H^2_0 \mstar_2  H^1_2)\mstar^1 Y)\star X\\&=
(-1)^{|X||Y|}((H^2_0 \mstar_1 H^1_2)\mstar^1 Y)\mstar X=
(-1)^{|X||Y|}Y\boxdot X.
\end{align*}
Finally, the proof of associativity is pictorially illustrated on Figure \ref{fig:calc}.

The products $\boxdot$ on $\bM^\bal$ and $\diamond$ on $\lb\bM^\ast\rb^\bal$ are related via the chain map $\lb \cH^2_0\rb^\carrl\colon \bM^{\bal}\to  \lb\bM^\ast\rb^\bal$ in the following way.
\begin{prp} 
The map induced by $\lb \cH^2_0\rb^\carrl$ on homology is a homomorphism of rings $\Phi\colon\lb H(\bM^\bal),\boxdot\rb\to \lb H(\lb\bM^\ast\rb^\bal),\diamond\rb$. Moreover, if there is a unit $e$ for $\boxdot$ on $H(\bM^\bal)$, i.e.~$e\boxdot X = X\boxdot e$, for all $X\in H(\bM^{\bal})$, then $\Phi$ is injective.
\end{prp}
Indeed, if $X,Y\in \bM^\bal$ then, taking into account that $|\Phi(X)|=|X|-1$, we calculate
\begin{align*}
\Phi(X)\diamond\Phi(Y)&=(-1)^{|X|+1}(H^2_0\star Y)\mstar((H^2_0\star X)\mstar_1 H^1_2)\\
&=H^2_0\star (((Y\star H^2_0) \mstar_1 H^1_2)\star X)\\
&=(-1)^{|X||Y|}H^{2}_0\mstar(((H^{2}_0\mstar_1H^{1}_{2})\mstar^{1}Y)\mstar X)\\
&=(-1)^{|X||Y|}\Phi(Y\boxdot X)=\Phi(X\boxdot Y).
\end{align*}
To verify the injectivity of $\Phi$ we use one of equivalent definitions of $\boxdot$ from Theorem \ref{thm:product-properties}, part 1:
\[
X=X\boxdot e=(-1)^{|X|+1}\lb \cH^1_2\rb^{\dua}\lb\Phi(X),e\rb,
\] 
and hence, $\Phi(X)=0$ implies $X=0.$

Using the isomorphism $\ul\beta\colon M^\cyc\to \bM^\bal$, $\ul\beta(X)=\Ee(X\sigma^{-1}\hbar^{-1})$ defined above we transport the product to $M^\cyc$. Given $X,Y\in M^\cyc$, we define
\begin{align}\notag
X\boxtimes Y&:=\ul\beta^{-1}\lb\ul\beta(X)\boxdot\ul\beta(Y)\rb\\\notag
&=
\lb  ((H^2_0\mstar\limits_1 H^1_2) \mstar^1\Ee(X\sigma^{-1}\hbar^{-1}))\star \Ee(Y\sigma^{-1}\hbar^{-1})\rb\hbar\sigma\\\label{eq:prodonM}
&=(-1)^{(n-2)|Y|+1} ((H^2_0\mstar_1 H^1_2) \mstar^1\Ee(X))\star \Ee(Y) \sigma^{-1}\hbar^{-1}.
\end{align}
 Note that the homomorphism $\ul\beta$ shifts the grading by $2-n$, and hence the commutation equation for $\boxdot$ translates to the equation $$X\boxtimes Y=(-1)^{(|X|-n)(|X|-n)}Y\boxtimes X.$$  Let us also point out that the degree of  the  product $\boxtimes$ is equal to $-n$.

\section{Legendrian algebra constructions}\label{sec:leg-algebra}
In this section we first reformulate the main result relating Legendrian homology algebra and symplectic homology \cite[Corollary 5.7]{BEE1} in the terminology of Section \ref{sec:alg-prelim}, and second give a brief sketch of the geometry that enters the construction of the full rational Hamiltonian in the framework of Legendrian SFT with focus on the parts of it that enter the expression of the product. 

\subsection{The complex $LH^{\Ho}(\Lambda)$ expressed as $M^{\cyc}$}\label{ssec:LHHo}
Let $(\oX_0,\lambda)$ be a Liouville domain with contact boundary $(Y_0,\ker(\lambda))$ and let $\Lambda=\Lambda_{1}\cup\dots \Lambda_{m}\subset Y_0$ be a union of Legendrian spheres. Let $\Cc$ be the set of Reeb chords of $\Lambda$ and consider the decomposition $\Cc=\bigcup_{1\le i,j\le m}\Cc_{ij}$, where a chord in $\Cc_{ij}$ starts on $\Lambda_{j}$ and ends on $\Lambda_{i}$.
Let $V$ be the graded module $K\la\Cc\ra$ generated by the set $\Cc$ with decomposition $V=\bigoplus_{1\le i,j\le m}V^{(ij)}$, where $V^{(ij)}$ is generated by $\Cc_{ij}$. Then the Legendrian homology algebra $LHA(\Lambda)$ is the tensor algebra $A$ associated to $V$ as in \eqref{eq:defA} and the differential $d_{LHA}$ gives a differential graded algebra $(A,d):=(LHA(\Lambda), d_{LHA})$, see \cite[Section 4.1]{BEE1}. Let $M=M(A)$ be the chain complex defined in \eqref{eq:defM}. The following result is a consequence of \cite[Remark 7.4]{BEE1}.

\begin{lma}\label{lm:LHHo} 
The chain complex $M^{\cyc}$ is canonically isomorphic to the chain complex $LH^{\Ho}(\Lambda)$ defined in \cite[Section 4.5]{BEE1}.
\end{lma}
In \cite[Corollary 5.7]{BEE1} it was shown that if $S\bbH(X_0)=0$ (e.g.~if $X_0$ is a subcritical Weinstein manifold) and if $X$ is obtained from $X_0$ by attaching Lagrangian $n$-handles along $\Lambda$, then  there is a quasi-isomorphism  
\begin{equation}\label{eq:qiso}
\Phi\colon SH(X)\to LH^{\Ho}(\Lambda)=M^{\cyc}\stackrel{\ul\beta}{=}\bM^{\bal},
\end{equation}
where $SH(X)$ denotes a chain complex with homology equal to $S\bbH(X)$.

\subsection{Legendrian rational SFT}
The differential in the Legendrian algebra $A=LHA(\Lambda)$ is defined by counting elements of moduli spaces of rigid holomorphic disks in $\R\times Y_0$, with boundary on $\R\times\Lambda$, and with exactly one positive and several negative punctures, see e.g.~\cite[Section 4.1]{BEE1}. The formalism of Legendrian rational SFT comes from algebraic structures associated with holomorphic disks with an arbitrary number of positive and negative punctures. As is well known, in order to construct such a theory one needs to handle boundary cusp-degenerations of  holomorphic disks. The problem of boundary cusp-degeneration can be solved by incorporating string topology operations into the SFT formalism. The corresponding theory has not yet been constructed in sufficient generality for it to be applicable in the setting of this paper. (See, however \cite{Ng} for the case of $1$-dimensional Legendrian knots in $\R^{3}$).
 
However, the algebraic formalism described in Section \ref{sec:alg-prelim} requires a different  version of  the full rational Hamiltonian $\cH$,   in the spirit of \cite{E:rsft}. It arises naturally when one considers moduli spaces of holomorphic disks with boundaries on multiple parallel copies of $\R\times\Lambda$. 

Let us first consider the term $\cH^2_0$. This is the only term besides $\cH^1_2$ which is needed for defining
the product on $M^\cyc$, see Figure \ref{fig:copies}. We recall here that $\cH^1_k$ for any $k\geq 1$ is determined by the differential on the Legendrian algebra $A$.

For each $i=1, \ldots, m$ fix a point $y_i \in \Lambda_i$. 
Let $Z=\llb \ul u_1w_1u_2w_2\rrb\in \lb\bU^2_0\rb^{\bal}$, where 
$w_j=c_j^1\dots c_j^{r_j}\in A$ and $u_j\in\wh\Cc^\ast \cup\Xx^\ast$, $j=1,2$, be an excited cyclic monomial. We associate to it a moduli space $\MM(Z)=\MM'(Z)\cup\MM''(Z)$. Here the moduli space $\MM'(Z)$ consists of holomorphic disks anchored in $X_0$, see \cite[Section 2.2]{BEE1},
$$
f\colon(D,\pa D \setminus\mathbf{z})\to (\R\times Y_0,\R\times \Lambda), 
$$
where
$$
\mathbf{z}=\left\{
z_{u_1}\,,\,z_{c_1^{1}},\dots,z_{c_1^{r_1}}\,,\, z_{u_2}\,,\, z_{c_2^{1}},\dots, z_{c_2^{r_2}}\right\}\subset\pa D
$$
is a set  of cyclically ordered boundary punctures on $\pa D$, with the following properties:
 \begin{itemize}
\item[($\mathrm{ch}^{-}$)]  at the puncture $z_{c^i_j}$, $f$ is asymptotic to the chord $c^i_j$ at $-\infty$;
\item[($\mathrm{ch}^{+}$)] if $u_j=\wh c^{\,\ast}\in \wh \Cc^\ast$ then  at the puncture
$z_{u_j}$  $f$ is asymptotic  to the chord $c$ at $+\infty$;
\item[$(\mathrm{pt}^{+})$] if $u_j=x_i^\ast\in \Xx^\ast$ then, $f$ extends continuously over $z_{u_j}$, and $f(z_{u_j})\in \R\times \{y_i\}$.
 \end{itemize}
 
The moduli space $\MM''(Z)$ consists of holomorphic disks connected by Morse flow lines. To define it, choose auxiliary Morse functions $\phi_i\colon\Lambda_i\to\R$, $i=1,\dots,m$, each with exactly two critical points. For $(s,t)\in\{0,1,\dots,r_1\}\times\{0,1,\dots,r_2\}$ define the auxiliary  moduli space $\MM_{(s,t)}(Z)$ of pairs $(f_1,f_2)$ of holomorphic disks anchored in $X_0$,
$$
f_\sigma\colon(D,\pa D \setminus\mathbf{z}_\sigma)\to (\R\times Y_0,\R\times \Lambda), \;\;\sigma=1,2,
$$ 
where   
\begin{align*}
\mathbf{z}_1&=\left\{z_{u_1},\,z_{c_1^{1}},\dots,z_{c_1^{s}}\,,\,\zeta_1,
  z_{c_2^{t+1}},\dots, z_{c_2^{r_2}}\right\}\subset\pa D,\\
\mathbf{z}_2&=\left\{z_{u_2},\,z_{c_2^{1}},\dots,z_{c_2^{t}}\,,\,\zeta_2,
  z_{c_1^{s+1}},\dots, z_{c_2^{r_2}}\right\}\subset\pa D
\end{align*}
are sets of cyclically ordered boundary punctures on $\pa D$, with the following properties. At $z_{u_{\sigma}}$, $f_\sigma$ satisfies 
condition $(\mathrm{ch}^{+})$ if $u_\sigma\in\wh \Cc^{\ast}$ or condition $(\mathrm{pt}^{+})$ if $u_{\sigma}\in\Xx^{\ast}$. At $z_{c_i^{j}}$, $f_\sigma$ satisfies condition $(\mathrm{ch}^{-})$. At the puncture $\zeta_\sigma$, $\sigma=1,2$,  $f_\sigma$ extends continuously and the following condition holds:
\begin{itemize}
\item[$(\mathrm{mo})$] The points $f_1(\zeta_1)$ and $f_2(\zeta_2)$ lie in the same component $\Lambda_i\subset\Lambda$ and there exists a trajectory of $-\nabla \phi_i$ which starts at $f_1(\zeta_1)$ and ends at $f_2(\zeta_2)$, where $\phi_i\colon \Lambda_i\to\R$ is the auxiliary Morse function on $\Lambda_i$. Here the gradient is taken with respect to  the metric induced by the symplectic form and the almost complex structure. 
\end{itemize}
Finally, define   
$$
\MM''(Z):=\bigcup_{(s,t)\in\{0,\dots,r_1\}\times\{0,\dots,r_2\}}\MM_{(s,t)}(Z).
$$

%Note also that we use the information about the location of the excited puncture at this point to determine. Role of $y_i$.

Given $u\in \wh \Cc\cup\wh \Cc^{\ast}$, let $\wt u=a$ if $u=\wh a$ or $u=\wh a^{\,\ast}$. Define $|\wt x_j|=|\wt x_j^{\ast}|=-1$, $j=1,\dots,m$. Then we have the dimension formula  
\[
\dim(\M(Z))= |\wt u_1|+|\wt u_2|-|w_1|-|w_2|-(n-3).
\]
In particular in terms of the grading on $U$, for the balanced monomial $\sigma^{2}\hbar^{-1}Z$ we have
\[
\dim(\M(Z))=-|\sigma^{2}\hbar^{-1}Z|.
\]
Define $\cH^2_0\in \lb \bU^2_0\rb^{\bal}$ as $\Ee(h^2_0)$, where
\begin{equation}\label{eq:Legham}
h^2_0\;=\sigma^{2}\hbar^{-1}\sum\limits_{|\sigma^{2}\hbar^{-1} Z|=-1}\#\M(Z)\,Z,
\end{equation}
where the sum is taken over all balanced monomials $\sigma^{2}\hbar^{-1} Z$ which additively generate $\lb \bU^2_0\rb^\bal$ and for which $|\sigma^{2}\hbar^{-1} Z|=-1$.
Here $\#\M(Z)$ is the algebraic number of $1$-dimensional components of the moduli space $\M(Z)$. As will be briefly discussed below, the above definition of $H^{2}_0$ can generalized to more than two positive punctures and as a consequence we get the following result. 

\begin{prp}\label{prop:Legendrian-master}
There exists an element $\cH'' \in \bU^\bal$, with trivial component in  $\lb \bU^1\rb^\bal\oplus \lb \bU^2_0\rb^\bal$, such that $\cH'=\cH^{2}_0+\cH''$ complements $\cH^1$ to the full rational SFT Hamiltonian (i.e.~if $\cH=\cH^{1}+\cH^2_0+\cH''$ then $d\cH+\cH\star \cH=0$).
\end{prp}

\begin{figure}
\centering
\includegraphics[width=.5\linewidth]{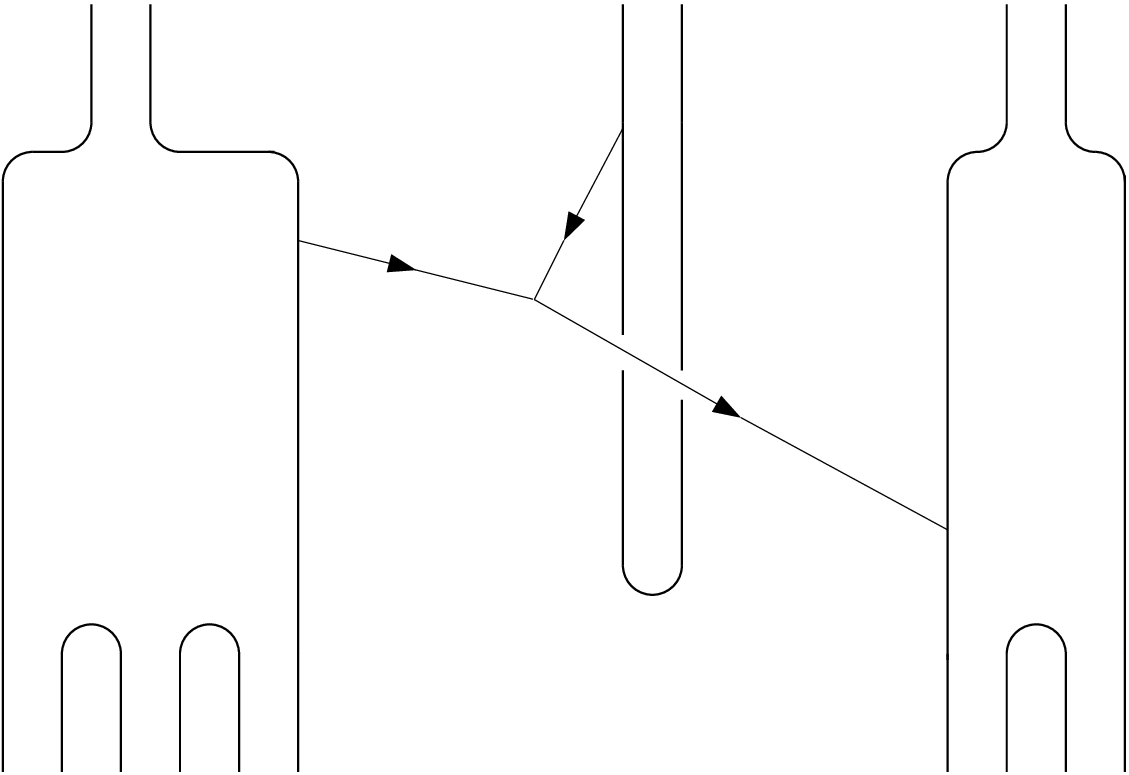}
\caption{A generalized holomorphic disk in $Y_0\times\R$ with boundary on $\R\times\Lambda$.}
\label{fig:gendisk}
\end{figure}

Here
the components $\cH^p_q\in\lb \bU^p_q\rb^{\bal}$, $p>2$, of the full rational Hamiltonian $H$ satisfies the relation $\cH^{p}_q=\frac1{q!}\Ss^{q}(h^{p}_0)$, $h^{p}_0=\frac1p \pi(H^{p}_0)$, and are defined through a construction which generalizes that used in the definition of $\cH^2_0$ above. These higher order components of the Hamiltonian are not needed for the definition of the product (but some of them are needed for proving its properties) and will not be defined in this paper. We give however a brief description: $\cH^p_q$ is defined using moduli spaces $\MM(Z)$ of so-called {\em generalized holomorphic disks}. Generalized holomorphic disks are combined objects consisting of holomorphic disks with boundary on one copy of $\Lambda$ connected by gradient flow trees associated to auxiliary Morse functions on $\Lambda$, an example is depicted in Figure \ref{fig:gendisk}, see \cite{EESa} for similar considerations. The elements of the moduli space $\MM''(Z)$ defined above are examples of generalized holomorphic disks where the gradient flow tree consists of only one gradient trajectory. A more general configuration is shown in Figure \ref{fig:gendisk}. More precisely, if $Z=\ul u_{1}w_1u_2w_2\dots u_k w_k$ is a monomial in $\bU^{p}_q(0,0)$, where $u_i\in\wh \Cc\cup\Xx\cup\wh \Cc^{\ast}\cup\Xx^{\ast}$ and where $w_i\in A$, $i=1,\dots,k$, then we consider the moduli space $\MM(Z)$ of generalized holomorphic disks with asymptotic data according to $Z$, i.e.~positive punctures at $\wt u_i$ if $u_i\in \wh\Cc^{\ast}\cup\Xx^{\ast}$, negative punctures at $\wt u_j$ if $u_j\in\wh\Cc\cup\Xx$, auxiliary negative punctures at the Reeb chords in the monomials $w_j$, and these punctures appearing along the boundary of the disk in the cyclic order given by $Z$. The dimension of $\MM(Z)$ is given by
\[
\dim(\MM(Z))=(n-3)+\sum_{u_i\in\wh\Cc^{\ast}\cup\Xx^{\ast}}\bigl(|\wt u_i|-(n-3)\bigr)-\sum_{u_j\in\wh\Cc\cup\Xx}|\wt u_j|-\sum_{j=1}^{m}|w_j|.
\] 
In particular for $\sigma^{p-q}\hbar^{-1}Z\in \lb \bU^{p}_q\rb^{\bal}$ we have
\[
\dim(\MM(Z))=-|\sigma^{p-q}\hbar^{-1}Z|
\] 
and $\cH^{p}_q$ is defined as
\[
\cH^{p}_{q}=\sigma^{p-q}\hbar^{-1}\sum_{|\sigma^{p-q}\hbar^{-1} Z|=-1}\#\M(Z)\,Z,
\]
where the sum ranges over all $\sigma^{p-q}\hbar^{-1} Z\in \lb \bU^{p}_{q}\rb^{\bal}$ with $|\sigma^{p-q}\hbar^{-1}Z|=-1$. 

\begin{rmk}
In this geometric context, the excitation operator $\Ee$ has the following meaning. Given a generalized holomorphic disk with boundary in the ordered union of (symplectizations of) parallel copies of a Legendrian manifold $\Lambda$ we always insist that the boundary arcs of the disk between punctures corresponding to generators from $\wh\Cc\cup\wh\Cc^\ast\cup\XX\cup\XX^\ast$ be mapped to different copies of $\Lambda$. We also require that when going around the boundary of the disk in the counter-clockwise direction we jump at every such puncture except one to a higher copy of $\Lambda$ in terms of the natural ordering corresponding to the direction of the Reeb vector field. The unique puncture where we jump down corresponds to the excited generator of the monomial corresponding to the boundary punctures. Here the moduli spaces corresponding to different assignment of the distinguished ``jump-down" punctures are canonically diffeomorphic, and the operator $\Ee$ just associates with this universal moduli space the moduli spaces corresponding to all possible assignments of the ``jump-down" vertex.
\end{rmk}

\begin{rmk}\label{r:gluingsign}
The geometric motivation for the ``gluing sign'' in the equations $\ux_j^\ast(x_j)=(-1)^{n-1}\hbar$, $j=1,\dots,m$ is as follows. The algebraic variables $x_j$ and $x^{\ast}_j$, $j=1,\dots,m$ has two dual geometric meanings as follows: when $x_j^{\ast}$ is glued to $\ux_j$, $x_j^{\ast}$ corresponds to a positive puncture at a short Reeb chord between copies of $\Lambda$ at the minimum of the Morse function $\phi\colon\Lambda\to\R$ and $\ux_j$ corresponds to a negative puncture at this chord, when $\ux_j^{\ast}$ is glued to $x_j$, $\ux_j^{\ast}$ corresponds to a negative puncture at the maximum of $\phi$ and $x_j$ to a positive puncture at this maximum. In a coherent orientation scheme there is a sign difference (a factor $(-1)^{n-1}$) between the gluings at these two types of punctures.    
\end{rmk}
    
\section{The product and the BV-operator via Legendrian algebra}\label{sec:BV-Leg}
Let $X_0$, $\Lambda\subset Y_0$, and $X$ be as above. Let $A=LHA(\Lambda)$ and $M=M(A)$. Then \eqref{eq:qiso} gives a quasi-isomorphism $\Phi\colon SH(X)\to M^{\cyc}\approx \bM^{\bal}$. In this section we present the operations on $M^{\cyc}$ which correspond to the product and the BV-operator on $SH(X)$ under $\Phi$.
\subsection{The product}
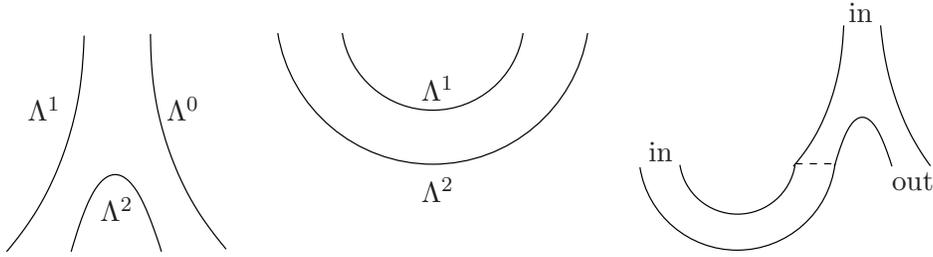
\begin{figure}
\labellist
\small\hair 2pt
\pinlabel $\Lambda^1$ at 31 114
\pinlabel $\Lambda^0$ at 140 114
\pinlabel $\Lambda^2$ at 88 35
\pinlabel $\Lambda^1$ at 340 130
\pinlabel $\Lambda^2$ at 340 50
\pinlabel $\text{in}$ at 672 190
\pinlabel $\text{out}$ at 713 57
\pinlabel $\text{in}$ at 515 80 
\endlabellist
\centering
\includegraphics[width=.9\linewidth]{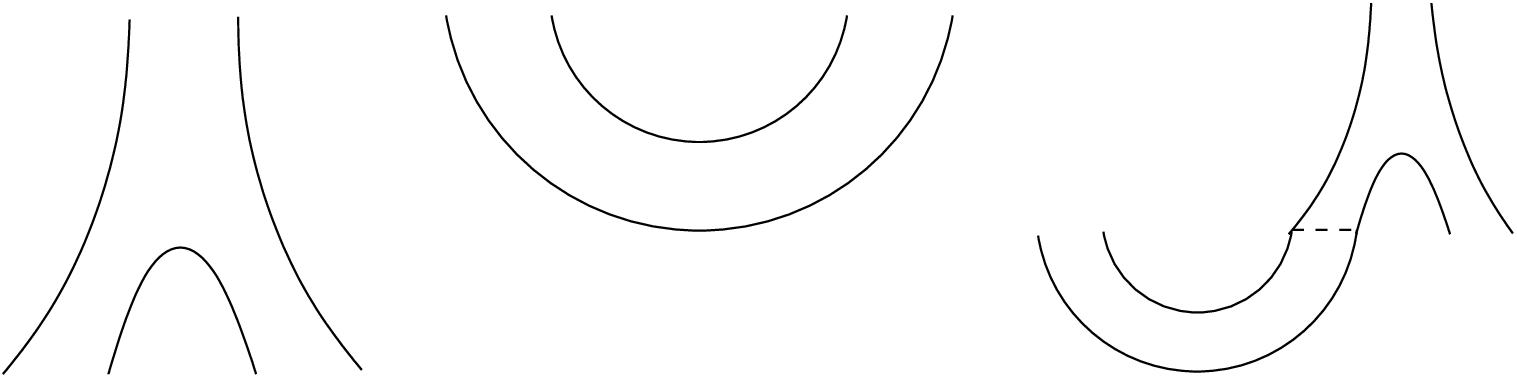}
\caption{Holomorphic disks used in the expression of the product, $\Lambda^j$, $j=0,1,2$ are parallel copies of $\Lambda$. The leftmost disk has three punctures at mixed (i.e.~connecting different components) Reeb chords, one positive and two negative, and auxiliary negative punctures at pure (i.e.~connecting a component to itself) chords which are not shown in the picture. The middle disk has two mixed positive punctures and auxiliary pure negative punctures, the latter again not shown. The rightmost picture shows how the disks are joined at a common mixed chord between $\Lambda^2$ and $\Lambda^1$ to form the tensor, with inputs and output as indicated, which represents the product.}
\label{fig:copies}
\end{figure}

Let $\cH\in \bU^{\bal}$ be the full rational Hamiltonian determined by $\Lambda\subset Y_0$, see \eqref{eq:Legham}. Then Proposition  \ref{prop:Legendrian-master} implies that $\cH^2_0$ satisfies \eqref{eq:H20}. Hence, according to Theorem \ref{thm:product-properties} and \eqref{eq:prodonM}, there exists a product  $\boxtimes\colon M^\cyc\otimes M^\cyc\to M^\cyc$ of degree $-n$ which is associative and commutative on homology and which is defined by conjugating  one of the following homotopically equivalent formulas for the product on $\bM^{\bal}$ with $\ul\beta$:
\begin{align}\notag
X\boxdot Y&\;=\;(\cH^2_0\mstar_1 \cH^1_2)^\downarrow(X,Y)\;\simeq\; -(\cH^2_0\mstar_2 \cH^1_2)^{\downarrow }(X,Y)\\\notag
&\simeq \;
(-1)^{|X|}\lb \cH^1_2\rb^{\dua}\circ\lb\lb \cH^2_0\rb^{\carrl}(X)\otimes Y\rb\\\label{eq:prodfmlas}
&\simeq\;(-1)^{|X||Y|+|X|+|Y|}\lb \cH^1_2\rb^{\uda}\circ\lb X\otimes\lb \cH^2_0\rb^{\carrl}(Y) \rb,
\end{align}
see Figure \ref{fig:copies}.
\begin{thm} \label{thm:SH-product}
If $S\bbH(X_0)=0$ (e.g.~if $X_0$ is subcritical) then,
under the map on homology induced by the quasi-isomorphism $\Phi\colon SH(X)\to M^\cyc$, see \eqref{eq:qiso}, the pair of pants product on $S\bbH(X)$ corresponds to the product $\boxtimes$ on $H(M^{\cyc})$.
\end{thm}

\subsection{The BV-operator}
Let $A=LHA(\Lambda)$ and consider the linear map $\Delta\colon M^{\cyc}\to M^{\cyc}$ of degree $1$ defined as follows on monomials $\llb uw\rrb$, $u\in \wh \Cc\cup\Xx$, $w\in A$:
\[
\Delta(\llb u w\rrb)= 
\begin{cases}
S(w) &\text{if } u\in \Xx,\\
0 &\text{if }u\in \wh \Cc. 	
\end{cases}
\]
\begin{prp}
The map $\Delta$ is a chain map and if $S\bbH(X_0)=0$ then,
under the map on homology induced by the quasi-isomorphism $\Phi\colon SH(X)\to M^\cyc$, see \eqref{eq:qiso}, the BV-operator on $S\bbH(X)$ corresponds to the operator $\Delta$ on $H(M^{\cyc})$.
\end{prp}

\section{Example: $SH(T^*S^n)$}\label{sec:examples}
In this section we compute the symplectic homology with product and BV-operator for cotangent bundles of spheres presented as the result of attaching a Lagrangian $n$-handle to the Legendrian unknot in the boundary of the ball. 

\begin{figure}
\labellist
\small\hair 2pt
\pinlabel $\Lambda$ at 125 218
\pinlabel $y$ at 595 123
\endlabellist
\centering
\includegraphics[width=.7\linewidth]{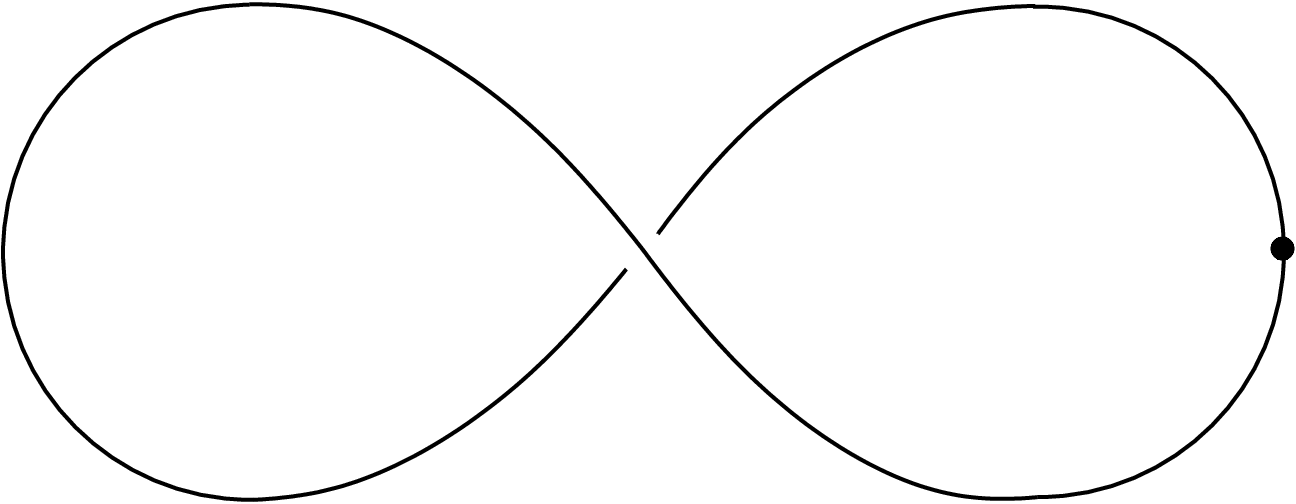}
\caption{The Legendrian unknot.}
\label{fig:unknot}
\end{figure}

Consider the Legendrian unknot in $\R^{2n-1}$ with coordinates $(q,p,z)\in\R^{n-1}\times\R^{n-1}\times\R$ and contact form $dz-p\cdot dq$. Figure \ref{fig:unknot} shows the projection of the Legendrian unknot $\Lambda$ for $n=2$ to the $qp$-plane.
To define the Legendrian unknot of dimension $n-1$ consider the $1$-dimensional unknot $\Lambda$ with Reeb chord over $q=p=0$ and symmetric with respect to the rotation by a half turn around the $z$ axis.  Include $\Lambda$ into the $q_1p_1z$-subspace of $\R^{2n-1}$ and take the $(n-1)$-dimensional unknot as the orbit of $\Lambda$ under the action by $SO(n-1)$ given by $\sigma(q,p,z)=(\sigma\cdot q,\sigma\cdot p,z)$ for $\sigma\in SO(n-1)$. For simplicity, we let $\Lambda$ denote the Legendrian unknot in $\R^{2n-1}$ for any $n$.  

Consider $\Lambda\subset\R^{2n-1}$ and include it into a Darboux chart in the boundary of the $2n$-ball $B^{2n}$. Then attaching a Lagrangian $n$-handle to $B^{2n}$ along $\Lambda$ gives a symplectic manifold which is symplectomorphic to the cotangent bundle $T^{\ast} S^{n}$ of the $n$-sphere. As explained in \cite{BEE1}, in order to determine $LHA(\Lambda)$ it suffices to study Reeb chords and holomorphic disks that are contained in the Darboux chart. The unknot $\Lambda$ has only one Reeb chord which we denote $a$ and $|a|=n-1$. 

The moduli space of holomorphic disks with positive puncture at $a$ and a marked point on the boundary evaluates with degree $\pm 1$ to $\Lambda$. To see this, note that for $n=2$, holomorphic disks that are rigid up to translation correspond to polygons with boundary on the knot diagram and with convex corners, where at positive punctures the boundary orientation goes from the lower to the upper strand and vice versa at negative punctures. There are thus two disks with evaluation maps that give a degree $\pm 1$ map since we use the null-cobordant spin structure on the circle to orient moduli spaces. For $n>2$ we note that we may take the almost complex structure to be invariant with respect to the $SO(n-1)$ action and it follows that the moduli space of holomorphic disks with positive puncture at $a$ is diffeomorphic to $S^{n-2}\times\R$ and that, when equipped with a marked point on the boundary,  it evaluates with degree $\pm1$.  

The above description of moduli spaces implies in particular that the differential $d$ on $A=LHA(\Lambda)=\Q\la a\ra$ is trivial, i.~e., $d=0$. Consequently, the differential on $\bM^{\bal}$ satisfies 
\[
d_{M} X=dX+[\cH_1^{1},X]=[\cH_1^{1},X], 
\]
where 
\begin{equation}\label{eq:unknotH11}
\cH_1^{1}=\lb \llb a \ux \hbar^{-1}\wh a^{\ast}\rrb - \llb \ux a \hbar^{-1}\wh a^{\ast}\rrb\rb.
\end{equation}

Thus $H(\bM^{\bal})$ and $H(M^{\cyc})$ (recall the chain isomorphism $\ul\beta\colon M^{\cyc}\to \bM^{\bal}$, see \eqref{eq:beta}) are generated by the cycles listed in Tables \ref{tab:SHn-1odd} and \ref{tab:SHn-1even} when $n-1$ is even and odd, respectively. 
\begin{table}[htp]
\begin{center}
{\renewcommand{\arraystretch}{1.5}
\renewcommand{\tabcolsep}{0.2cm}
\begin{tabular}{|c|c||c|c|}
\hline
$H(\bM^{\bal})$ & degree & $H(M^{\cyc})$ & degree\\
\hline
$\llb \ux\rrb \hbar^{-1}\sigma^{-1}$ & $-(n-2)$ & $\llb x\rrb$ & $0$\\
\hline
$\llb \ux a^{2k+1}\rrb\hbar^{-1}\sigma^{-1}$ & $2k(n-1)+1$ & $\llb xa^{2k+1}\rrb$ & $(2k+1)(n-1)$\\
\hline
$\llb \wh \ua a^{2k}\rrb \hbar^{-1}\sigma^{-1}$ & $2k(n-1)+2$ & 
$\llb \wh a a^{2k}\rrb$ & $(2k+1)(n-1)+1$\\
\hline
\end{tabular}}
\end{center}
\caption{Homology of $\bM^{\bal}$ and $M^{\cyc}$ for $n-1$ odd, $k=0,1,2,\dots$.}\label{tab:SHn-1odd}
\end{table}

\begin{table}[htp]
\begin{center}
{\renewcommand{\arraystretch}{1.5}
\renewcommand{\tabcolsep}{0.2cm}
\begin{tabular}{|c|c||c|c|}
\hline
$H(\bM^{\bal})$ & degree & $H(M^{\cyc})$ & degree\\
\hline
$\llb \ux\rrb \hbar^{-1}\sigma^{-1}$ & $-(n-2)$ & $\llb x\rrb$ & $0$\\
\hline
$\llb \ux a^{k}\rrb \hbar^{-1}\sigma^{-1}$ & $(k-1)(n-1)+1$ & $\llb xa^{k}\rrb$ & $k(n-1)$\\ 
\hline
$\llb \wh \ua a^{k-1}\rrb \hbar^{-1}\sigma^{-1}$ & $(k-1)(n-1)+2$ & $\llb \wh aa^{k-1}\rrb$ & $k(n-1)+1$\\
\hline
\end{tabular}}
\end{center}
\caption{Homology of $\bM^{\bal}$ and $M^{\cyc}$ for $n-1$ even, $k=1,2,\dots$.}\label{tab:SHn-1even}
\end{table}

In order to compute the product we need to determine $\cH^{1}_{2}$ and $\cH^{2}_{0}$. Consider first $\cH^{1}_{2}$ (see Definition \ref{dfn:hamiltonian1}): 
\begin{align}\notag
h^{1}_{2}&=\frac12\Ss h^{1}_1 + \sigma^{-1}\llb xx\hbar^{-1}x^{\ast}\rrb\\\label{eq:H^1_2unknot}
&=\sigma^{-1}\lb \llb\wh ax\hbar^{-1}\wh a^{\ast}\rrb+ \llb x\wh a\hbar^{-1}\wh a^{\ast}\rrb  + (-1)^{n-1} \llb xx\hbar^{-1}x^{\ast}\rrb\rb,
\end{align}
see Remark \ref{r:signs} for a discussion of the sign in this formula.  
\begin{rmk}
In our description of $h$ in terms of generalized holomorphic disks, the disks which give rise to  \eqref{eq:H^1_2unknot}	are degenerate. The first two terms comes from the trivial strip of $a$ with an unconstrained Morse puncture lying in one of the two boundary components. The last term comes from the $\R$-family of three-punctured constant disks lying in $\{y\}\times\R$ with one puncture constrained to lie in $\{y\}\times\R$ and remaining two punctures unconstrained.    
\end{rmk}
Next, consider  $h^{2}_0$. Since the evaluation map from the moduli space of holomorphic disks with positive puncture at $a$ into $\Lambda$ has degree $\pm 1$ there is algebraically $\pm 1$ disk which satisfies a point constraint at $y\in \Lambda$ and we conclude from the description of generalized holomorphic disks that
\begin{equation}\label{eq:H^2_0unknot}
h^{2}_0=\sigma^{2}\llb \wh a^{\ast}x^{\ast}\rrb\hbar^{-1},
\end{equation}
see Remark \ref{r:signs} for a discussion of the sign in this formula. 

\begin{rmk}\label{r:signs}
Orientations of moduli spaces are fixed by orienting capping operators at Reeb chords, see \cite{EESori} for details. The overall sign in \eqref{eq:H^1_2unknot} is independent of the choice of capping operator for $a$ since variables associated to $a$ appears quadratically in $H^{1}_2$, the overall signs in \eqref{eq:H^1_2unknot} changes with the orientation data at $x$, and the overall sign in \eqref{eq:H^2_0unknot} changes, independently, with the orientations of capping operators at both $x$ and $a$. It follows that we can choose any sign combination in the two equations and in particular the one under consideration. 
\end{rmk}

To compute $H^{2}_0\mstar_1 H^{1}_{2}$ acting on $\bM^{\bal}\otimes\bM^{\bal}$ we use the components of $H^{1}_{2}$ in which the puncture at the second leg is excited. We will denote the sum of these by $\dot H^{1}_{2}$. The relevant parts of the Hamiltonian are thus, see \eqref{eq:H^2_0unknot} and \eqref{eq:H^1_2unknot}, (recall gradings: $|\wh a|=n$, $|\wh a^{\ast}|=n-3-n=-3$, $|x|=0$, and $|x^{\ast}|=n-3$)
\begin{align*}
H^{2}_0&=\sigma^{2}\lb\llb \wh a^{\ast}\hbar^{-1}
\ux^{\ast}\rrb-\llb x^{\ast}\hbar^{-1}\wh \ua^{\ast}\rrb\rb,\\
\dot H^{1}_{2}&=\sigma^{-1}\lb \llb\wh a\ux\hbar^{-1}\wh a^{\ast}\rrb+ \llb x\wh \ua\hbar^{-1}\wh a^{\ast}\rrb  +  (-1)^{n-1}\llb x\ux\hbar^{-1}x^{\ast}\rrb\rb.
\end{align*}
and we get
\begin{align*}
\cH^{2}_0\mstar_1 \dot\cH_{2}^{1}&= 
-H^{2}_0\;\mstar_1\;\lb\llb \wh a\ux\,\hbar^{-1}\wh a^{\ast}\rrb+\llb x\wh \ua\hbar^{-1}\wh a^{\ast}\rrb + (-1)^{n-1}\llb x\ux\hbar^{-1}x^{\ast}\rrb\rb\sigma^{-1}\\
&=
\sigma^{2}\lb\llb x^{\ast}\ux\hbar^{-1}\wh a^{\ast}\rrb
+(-1)^{n}\llb\wh a^{\ast}\wh \ua\hbar^{-1}\wh a^{\ast}\rrb
-\llb\wh a^{\ast}\ux\hbar^{-1}x^{\ast}\rrb
\rb\sigma^{-1}\\
&=
\sigma\lb-\llb x^{\ast}\ux\wh a^{\ast}\hbar^{-1}\rrb
+(-1)^{n-1}\llb\wh a^{\ast}\wh\ua\wh a^{\ast}\hbar^{-1}\rrb
+\llb\wh a^{\ast}\ux x^{\ast}\hbar^{-1}\rrb
\rb\\
&=
\sigma\lb-\llb \ux\wh a^{\ast}\hbar^{-1}x^{\ast}\rrb
+(-1)^{n-1}\llb\wh \ua\wh a^{\ast}\hbar^{-1}\wh a^{\ast}\rrb
+\llb \ux x^{\ast}\hbar^{-1}\wh a^{\ast}\rrb
\rb.
\end{align*}

Using \eqref{eq:prodfmlas} for the tensor of the product $\boxdot \colon H(\bM^{\bal})\otimes H(\bM^{\bal})\to H(\bM^{\bal})$ then gives the following expression for the product, with $\theta=\hbar^{-1}\sigma^{-1}$: If $n-1$ is odd then
\begin{align*}
&\llb \ux\rrb\theta \boxdot \llb\ux\rrb\theta=0,\quad
\llb \ux\rrb\theta \boxdot \llb \ux a^{2k+1}\rrb\theta=0,\quad
\llb \ux a^{2k+1}\rrb\theta\boxdot\llb \ux a^{2l+1}\rrb\theta=0,\\
&\llb \wh \ua a^{2k}\rrb\theta\boxdot\llb \wh \ua a^{2l}\rrb\theta=\llb \wh \ua a^{2(l+k)}\rrb\theta,\\
&\llb \wh \ua\rrb\theta \boxdot \llb \ux\rrb\theta=
\llb \ux\rrb\theta,\\
&\llb \wh \ua a^{2k}\rrb\theta \boxdot \llb \ux\rrb\theta=
\llb \ux a^{2k}\rrb\theta=
\tfrac12 d_{\tau M}\llb \wh \ua a^{2k}\rrb\theta\simeq 0 \text{ if $k>0$},\\
&\llb \wh \ua a^{2k}\rrb\theta \boxdot \llb \ux a^{2l+1}\rrb\theta=
\llb \ux a^{2(k+l)+1}\rrb\theta.
\end{align*}
and if $n-1$ is even then for $k,l\ge 0$
\begin{align*}
&\llb \ux a^{k}\rrb\theta\boxdot \llb \ux a^{l}\rrb\theta=0,\\
&\llb \wh \ua a^{k}\rrb\theta \boxdot \llb \ux a^{l}\rrb\theta= 
\llb \ux a^{k+l}\rrb\theta,\\
&\llb \wh \ua a^{k}\rrb\theta \boxdot\llb \wh \ua a^{l}\rrb\theta= 
\llb \wh \ua a^{k+l}\rrb\theta.
\end{align*}

We conclude from this calculation that the commutative product $\boxtimes\colon M^{\cyc}\otimes M^{\cyc}\to M^{\cyc}$ which corresponds to the pair of pants product on $S\bbH(T^{\ast}S^{n})$ has unit represented by the class $\llb\wh a\rrb\in M^{\cyc}$ and multiplication according to the above formulas after substituting $\theta$ with $1$ everywhere.

To finish the example we determine also the operator $\Delta\colon H(M^{\cyc})\to H(M^{\cyc})$ which corresponds to the $\BV$-operator on $S\bbH(T^{\ast}S^{n})$. If $n-1$ is odd then
\[
\Delta(\llb x\rrb)=0,\quad\Delta(\llb\wh a a^{2k}\rrb)=0,\; k\ge0,\quad
\Delta(\llb x a^{2k+1}\rrb)=\llb \wh a a^{2k}\rrb,\; k\ge 0,
\]
and if $n-1$ is even then
\[
\Delta(\llb x\rrb)=0,\quad\Delta(\llb \wh a a^{k-1}\rrb )=0,\;k\ge 1,\quad
\Delta(\llb x a^{k}\rrb)=\llb \wh a a^{k-1}\rrb,\;k\ge 1.
\]

\end{document}